\documentclass[11pt,twoside]{article}
\usepackage{amssymb,amsbsy,amsmath,amsfonts,amssymb,amscd,graphicx,color}
\usepackage[english,francais]{babel}
\newtheorem{theorem}{Theorem}[section]
\newtheorem{lemma}[theorem]{Lemma}
\newtheorem{proposition}[theorem]{Proposition}
\newtheorem{corollary}[theorem]{Corollary}
\newtheorem{definition}[theorem]{Definition\rm}
\newtheorem{remark}[theorem]{Remark}

\newcommand{\proof}{\noindent {\sl Proof.\/} }
\newcommand{\finproof}{\unskip\null\hfill$\square$\vskip 0.3cm}
\newcommand{\RR}{{\mathbb R}}
\newcommand{\R}{{\mathbb R}}
\newcommand{\LL}{{\rm L}}
\newcommand{\EE}{{\mathcal K}}

\newcommand{\WW}{{\mathcal W}}
\newcommand{\NN}{{\mathbb N}}
\newcommand{\CC}{{\mathcal C}}
\newcommand{\eps}{\varepsilon}
\newcommand{\FF}{\mathcal F}

\newcommand{\GG}{\mathcal G}

\newcommand{\ZZ}{\mathcal Z}

\newcommand{\un}{{\rm 1\kern -2.5pt l}}

\newcommand{\YY}{{\mathcal Y}}

\newcommand{\ee}{{\rm e}}
\newcommand{\II}{{\mathcal I}}
\newcommand{\dx}{\;{\rm d}x}
\newcommand{\dy}{\;{\rm d}y}
\newcommand{\dd}{\;{\rm d}}
\newcommand{\ds}{\;{\rm d}s}

\newcommand{\dt}{\;{\rm d}t}
\newcommand{\dr}{\;{\rm d}r}
\newcommand{\wto}{\rightharpoonup}
\newcommand{\email}[1]{{\small E-mail: {\textsf {#1}}}}
\newcommand{\http}[1]{{\small Internet: {\textsf {#1}}}}

\begin{document}\selectlanguage{english}
\title{Critical mass for a Patlak-Keller-Segel model \\  with degenerate diffusion in higher dimensions}

\author{Adrien Blanchet\footnote{EPI SIMPAF -- INRIA Futurs, Parc Scientifique de la haute Borne, F--59650 Villeneuve d'Ascq, France \& Laboratoire Paul Painlev\'e -- Universit\'e de Lille 1, F--59655 Villeneuve d'Ascq C\'edex, France.
\email{blanchet@ceremade.dauphine.fr},
\http{http://www.ceremade.dauphine.fr/$\sim$blanchet/}},\kern8pt
Jos\'e A. Carrillo\footnote{ICREA (Instituci\'o Catalana de
Recerca i Estudis Avan\c cats) and Departament de
Mate\-m\`a\-tiques, Universitat Aut\`onoma de Barcelona, E--08193
Bellaterra, Spain. \email{carrillo@mat.uab.es},
\http{http://kinetic.mat.uab.es/$\sim$carrillo/}}
\kern8pt\& \kern8pt Philippe Lauren\c{c}ot \footnote{Institut de
Math\'ematiques de Toulouse, CNRS UMR~5219 \& Universit\'e de
Toulouse, 118 route de Narbonne, F--31062 Toulouse C\'edex 9,
France. \email{laurenco@mip.ups-tlse.fr},
\http{http://www.mip.ups-tlse.fr/$\sim$laurenco/}}}
\date{\today}
\maketitle
\begin{abstract}
This paper is devoted to the analysis of non-negative solutions for a
generalisation of the classical parabolic-elliptic Patlak-Keller-Segel
system with $d\ge3$ and porous medium-like non-linear diffusion. Here,
the non-linear diffusion is chosen in such a way that its scaling and
the one of the Poisson term coincide. We exhibit that the qualitative
behaviour of solutions is decided by the initial mass of the
system. Actually, there is a sharp critical mass $M_c$ such that if $M
\in (0,M_c]$ solutions exist globally in time, whereas there are
blowing-up solutions otherwise. We also show the existence of
self-similar solutions for $M \in (0,M_c)$. While characterising the
eventual infinite time blowing-up profile for $M=M_c$, we observe that
the long time asymptotics are much more complicated than in the
classical Patlak-Keller-Segel system in dimension two. 
\end{abstract}
\section{Introduction}

In this work, we analyse qualitative properties of non-negative solutions for
the Patlak-Keller-Segel system in dimension $d\geq 3$ with
homogeneous non-linear diffusion given by
\begin{equation}\label{eq:sp}
\left\lbrace
\begin{array}{rll}
\displaystyle \frac{\partial u}{\partial t}(t,x)&={\rm div } \left[
\nabla u^m(t,x)-u(t,x)\nabla \phi(t,x)\right]\qquad &
t>0\,,\;x\in\RR^d\;,\vspace{.3cm}\\
\displaystyle -\Delta \phi(t,x)&=u(t,x)\;,\qquad
&t>0\,,\;x\in\RR^d\;,\vspace{.3cm}\\
u(0,x)&=u_0(x)\qquad &x\in\RR^d\, .
\end{array}
\right.
\end{equation}
Initial data will be assumed throughout this paper to verify
\begin{equation}\label{Eqn:Assumptions}
u_0\in \LL^1(\RR^d ; (1+|x|^2)\; \dd x)\,\cap\, \LL^\infty(\RR^d),
\quad\nabla u_0^m \in \LL^2(\R^d) \;\;\mbox{ and }\;\; u_0\ge0 \;.
\end{equation}
A fundamental property of the solutions to \eqref{eq:sp} is the formal
conservation of the total mass of the system 
\begin{equation*}
M:=\int_{\RR^d}u_0(x)\; \dd x=\int_{\RR^d}u(t,x)\; \dd x\quad
\mbox{for $t\ge 0$}\;. 
\end{equation*}
As the solution to the Poisson equation $-\Delta \phi=u$ is given up to
an harmonic function, we choose the one given by
$\phi=\EE\ast u$ with
\begin{equation*}
\EE(x)=c_d \;\frac{1}{|x|^{d-2}} \quad\mbox{and}\quad
  c_d:=\frac {1}{(d-2) \sigma_d}
\end{equation*}
where $\sigma_d:=2\,\pi^{d/2}/\Gamma(d/2)$ is the surface area of the
sphere $\mathbb{S}^{d-1}$ in $\RR^d$. This system has been
proposed as a model for chemotaxis-driven cell movement or in the
study of large ensemble of gravitationally interacting particles,
see \cite{horst,ChS04,BDP} and the literature therein.

We will concentrate on a particular choice of the non-linear
diffusion exponent $m$ in any dimension characterised for
producing an exact balance in the scaling of diffusion and
potential drift in equation \eqref{eq:sp}. To this end we use the
by-now classical scaling leading to the nonlinear Fokker-Planck
equation for porous media as in \cite{CT00}, that is, let us
define $\rho$ by $ \rho(s,y): = \ee^{dt}
u\left(\beta(t),\ee^{t}x\right)$ and $c:=\EE
\ast \rho$ with $\beta$ strictly increasing to be chosen. Then,
it is straightforward to check that
\begin{equation*}
\left\lbrace
\begin{array}{rll}
\displaystyle \frac{\partial \rho}{\partial s}(s,y)&={\rm div}
\left[ y \rho(s,y) + \beta'(t)\left\{ \ee^{-\lambda t}\nabla
\rho^m(s,y)- \ee^{-d t}\rho(s,y)\nabla c(s,y)\right\}\right]\quad &
s>0\,,\;y\in\RR^d\;,\vspace{.3cm}\\
\displaystyle -\Delta c(s,y)&=\rho(s,y)\;,\quad
&s>0\,,\;y\in\RR^d\;,\vspace{.3cm}\\
\rho(0,y)&=u_0(y)\ge 0\quad &y\in\RR^d\,,
\end{array}
\right.
\end{equation*}
with $\lambda=d(m-1)+2$. From this scaling, the only possible
choice of $m$ leading to a compensation effect between diffusion
and concentration is given by $\lambda=d$ or equivalently
\begin{equation*}
  m_d:=\frac{2\,(d-1)}{d}\;.
\end{equation*}
In that case, $\beta'(t)=\ee^{d t}$ determines the change of
variables and the final scaled equation reads:
\begin{equation}\label{eq:spscaled}
\left\lbrace
\begin{array}{rll}
\displaystyle \frac{\partial \rho}{\partial s}(s,y)&={\rm div }\left[
y \rho(s,y) + \nabla \rho^{m_d}(s,y)-\rho(s,y)\nabla
c(s,y)\right]\qquad & s>0\,,\;y\in\RR^d\;,\vspace{.3cm}\\
\displaystyle -\Delta c(s,y)&=\rho(s,y)\;,\qquad
&s>0\,,\;y\in\RR^d\;,\vspace{.3cm}\\
\rho(0,y)&=u_0(y)\geq 0\qquad &y\in\RR^d\,,
\end{array}
\right.\end{equation}
Note that the case $d=2$ and $m_2=1$ corresponds to the
Patlak-Keller-Segel system or to the classical
Smoluchowski-Poisson system in two dimensions with linear
diffusion~\cite{patlak,Keller-Segel-70}. In this case, a simple
dichotomy result have been shown 
in \cite{DP,BDP} improving over previous results in
\cite{Jager92,Nagai01}, namely, the behaviour of the solutions is
just determined by the initial mass of the system. More precisely,
there exists a critical value of the mass $M_c:=8\pi$ such that if
$0<M<M_c$ (sub-critical case) the solutions exist globally and if
$M>M_c$ (super-critical case) the solutions blow up in finite time.
Moreover, in the sub-critical case solutions behave self-similarly
as $t\to\infty$ \cite{bkln,BDP}. Finally, the critical case $M=M_c$
was studied in 
\cite{BCM} showing that solutions exist globally and blow up as a
Dirac mass at the centre of mass as $t\to\infty$. Solutions have
to be understood as free energy solutions, concept that we will
specify below.\medskip

In this work, we will show that a similar situation to the
classical PKS system in $d=2$,
although with some important differences, happens for the critical
variant of the PKS model in any dimension $d\geq 3$ reading as:
\begin{equation}\label{eq:spf}
\left\lbrace
\begin{array}{rll}
\displaystyle \frac{\partial u}{\partial t}(t,x)&={\rm div }
\left[ \nabla u^{m_d}(t,x)-u(t,x)
\nabla (\EE \ast u)(t,x)\right]\qquad & t>0\,,\;x\in\RR^d\;,\vspace{.3cm}\\
u(0,x)&=u_0(x)\geq 0\qquad &x\in\RR^d\, .
\end{array}
\right.
\end{equation}
We will simply denote by $m$ the critical exponent
$$m:=m_d=\frac{2\,(d-1)}{d}\in (1,2)\,,$$ as long as $d\geq 3$, in the rest of
the paper for notational convenience. The main tool for the
analysis of this equation is the following {\it free energy}
functional:
\begin{align*}
 t\mapsto \FF[u(t)]:&= \int_{\RR^d} \frac{u^m(t,x)}{m-1} -
 \frac12 \iint_{\RR^d \times \RR^d} \EE(x-y)\,u(t,x)\,u(t,y)\dd x \dd y\;\\
&= \int_{\RR^d} \frac{u^m(t,x)}{m-1} -
 \frac{c_d}{2} \iint_{\RR^d \times \RR^d}
\frac{1}{|x-y|^{d-2}}\,u(t,x)\,u(t,y)\dd x \dd y\;\nonumber
\end{align*}
which is related to its time derivative, {\it the Fisher
information}, in the following way: given a smooth positive
fast-decaying solution to~\eqref{eq:spf}, then
\begin{equation}
\label{Lemma:FreeEnergy}
\frac\dd{\dt}\FF[u(t)]=-\int_{\RR^d}u(t,x)\left|\nabla\left(\frac{m}{m-1}
u^{m-1}(t,x) -\phi(t,x) \right)\right|^2\dx \; .
\end{equation}
We will give a precise sense to this entropy/entropy-dissipation
relation below.

The system \eqref{eq:spf} can formally be considered a particular
instance of the general family of PDEs studied in
\cite{Carrillo-McCann-Villani03,
AmbrosioGigliSavare02p, Carrillo-McCann-Villani06}. The free
energy functional $\FF$ structurally belongs to the general class
of free energies for interacting particles introduced in
\cite{McCann97,
Carrillo-McCann-Villani03,Carrillo-McCann-Villani06}. The
functionals treated in those references are of the general form:
\begin{equation*}
{\cal E}[n]:=\int_{\RR^d} U[n(x)]\dx + \int_{\RR^d} n(x)\, V(x)\dx
+ \frac12 \iint_{\RR^d\times\RR^d} W(x-y)\, n(x)\,n(y)\dx \dy\;
\end{equation*}
under the basic assumptions $U:\R^+\to \R$ is a density of
internal energy, $V:\R^d\to\R$ is a convex smooth confinement
potential and $W:\R^d\to\R$ is a symmetric convex smooth
interaction potential. The internal energy $U$ should satisfy the
following dilation condition, introduced in McCann~\cite{McCann97}
\begin{equation*}
\lambda \longmapsto \lambda^d U(\lambda^{-d}) \qquad \text{is
convex non-increasing on $\R^+$}.
\end{equation*}
In our case, the interaction
potential is singular and the key tool of displacement convexity
of the functional fails, making the theory in the previous
references not useful for our purposes. Nevertheless, the free
energy functional plays a central role for this problem as we
shall see below. Before proceeding further, let us state the
notion of solutions we will deal with in the rest:

\begin{definition}[Weak and free energy solution]\label{def:wfes}
Let $u_0$ be an initial condition satisfying
\eqref{Eqn:Assumptions} and $T\in (0,\infty]$.

\begin{itemize}
\item[(i)] A \emph{weak solution} to \eqref{eq:spf} on $[0,T)$ with initial
condition $u_0$ is a non-negative function $u\in
C([0,T);\LL^1(\RR^d))$ such that $u\in \LL^\infty((0,t)\times\RR^d))$, $u^m\in
\LL^2(0,t;H^1(\RR^d))$ for each $t\in [0,T)$ and
\begin{multline}\label{eq:weaksol}
\int_{\RR^d} u_0(x)\,\psi(0,x)\dd x  =\int_0^T \!\!\int_{\RR^d} \left[ \nabla u^m(t,x) - u(t,x)\nabla\phi(t,x)\right]\cdot\nabla \psi(t,x)  \dd x \dd t \\
  -\int_0^T \!\!\int_{\RR^d }u(t,x)\,\partial_t\psi(t,x) \dd x \dd t
\end{multline}
for any test function $\psi \in \mathcal D([0,T)\times\RR^d)$ with
$\phi=\EE\ast u$. \item[(ii)] A \emph{free energy solution}
to~\eqref{eq:spf} on $[0,T)$ with initial condition $u_0$ is a
weak solution to~\eqref{eq:spf} on $[0,T)$  with initial condition
$u_0$ satisfying additionally: $u^{(2m-1)/2}\in
\LL^2(0,t;H^1(\RR^d))$ and
\begin{equation}\label{freeeneinq}
  \FF[u(t)] + \int_0^t \!\! \int_{\RR^d} \left|\left(\frac{2m}{2m-1}
\nabla u^{(2m-1)/2}(s,x) -u^{1/2}(s,x)\nabla \phi (s,x)\right)\right|^2\dd x\dd s \le \FF[u_0]
\end{equation}
for all $t\in(0,T)$ with $\phi=\EE\ast u$.
\end{itemize}
\end{definition}
In~\eqref{freeeneinq}, we cannot write the Fisher information factorised by $u$ as in~\eqref{Lemma:FreeEnergy} because of the lack of regularity of $u$. We note that both~\eqref{eq:weaksol} and \eqref{freeeneinq} are
meaningful. Indeed, the regularity required for $u$ implies that
the solution $\phi=\EE\ast u$ to the Poisson equation satisfies
$\phi\in\LL^\infty(0,t;H^1(\RR^d))$ for all $t\in (0,T)$. In
addition, it follows from~\eqref{eq:weaksol} by classical
approximation arguments that
\begin{equation}\label{evian}
\|u(t)\|_1 = \int_{\RR^d} u(t,x) \dd x =
\int_{\RR^d} u_0(x) \dd x = \|u_0\|_1 = M \;\;\mbox{ for }\;\;
t\in [0,T)\;.
\end{equation}

Let us point out that the existence of free energy solutions for a
related problem was essentially obtained in~\cite{S1,S2,O} where
the Poisson equation is replaced by $-\Delta \phi=u-\phi$. There,
the authors also show that the mass is the suitable quantity
for~\eqref{eq:spf} allowing for a dichotomy. Precisely, the author
shows that there exist two masses $0<M_1<M_2$ such that if
$0<M<M_1$ the solutions exist globally in time, while for $M>M_2$
there are solutions blowing up in finite time. The values of these
masses, are related to the sharp constants of the Sobolev
inequality.

Here, we will make a fundamental use of a {\it variant to the
Hardy-Littlewood-Sobolev (VHLS) inequality}, see
Lemma~\ref{lem:vhls}: for all $h\in \LL^1(\RR^d)\cap \LL^m(\R^d)$,
there exists an optimal constant $C_*$ such that
\begin{equation*}
\left|\iint_{\RR^d\times\RR^d}\frac {h(x)\,h(y)}{|x-y|^{d-2}} \dd
x\dd y\right|\le C_*\,\|h\|^m_{m}\,
\|h\|^{{2}/{d}}_{^1}\;.
\end{equation*}
This inequality will play the same role as the
logarithmic HLS inequality proved in \cite{CL} for the classical
PKS system in $d=2$ \cite{DP,BDP,BCM}. The VHLS inequality and the
identification of the equality cases allow us to give the first
main result of this work, namely, the following sharp critical
mass
\begin{equation*}
  M_c := \left[\frac{2}{(m-1)\,C_{*}\,c_d} \right]^{d/2}\;
\end{equation*}
for equation \eqref{eq:spf}. More precisely, we will show that
free energy solutions {\sl exist globally for $M\in (0, M_c]$} while
there are finite time blowing-up solutions otherwise. However, the
long time asymptotics of the solutions is much more complicated
compared to the classical PKS system in two dimensions. The main
results of this work and the open problems related to large times
asymptotics can be summarised as follows:
\begin{itemize}
\item Sub-critical case: $0<M<M_c$, solutions exist globally in
time and there exists a radially symmetric compactly supported
self-similar solution, although we are not able to show that it
attracts all global solutions. See
Proposition~\ref{prop:existenceglobalsol},
Theorem~\ref{thm:rescaledminimiser} and Corollary~\ref{cor:last}.

\item Critical case: $M=M_c$, solutions exist globally in time,
see Proposition~\ref{prop:gloexistcritic}. There are infinitely
many compactly supported stationary solutions. The
second moment of solutions is non-decreasing in time, with two
possibilities we cannot exclude: either is uniformly bounded in
time or diverges. Moreover, the $\LL^m$-norm of the solution
could be divergent as $t\to\infty$ or a diverging sequence of
times could exist with bounded $\LL^m$-norm. However, we show a
striking difference with respect to the classical PKS system in
two dimensions~\cite{BCM}, namely, the existence of global in time
solutions not blowing-up in infinite time. We will comment further
on these issues in Section~\ref{sec:does}.

\item Super-critical case: $M>M_c$, we prove that there exist solutions, corresponding to initial data with negative free energy, blowing up in finite time, see Proposition~\ref{lem:superfinite}. However, we cannot exclude the possibility that solutions with positive free energy may be global in time.
\end{itemize}

The results are organised as follows. Section~\ref{sec:existence} shows a key maximal time of existence criterion for free energy solutions of equation~\eqref{eq:spf}. This criterion improves over the results in~\cite{S1,S2} since it is only based on the
boundedness or unboundedness in time of the $\LL^m$-norm of the solutions and
it has to be compared to a similar criterion based on the
logarithmic entropy in the classical PKS system in two dimensions
obtained in~\cite{BCM}. Section 3 is devoted to the variational
study of the minimisation of the free energy functional over the
set of densities with a fixed mass. With that aim the proof of the
VHLS inequality and the identification of the
equality cases are performed. Section 4 uses this variational
information to show the above main results concerning the
dichotomy, the global existence for $M< M_c$ and the
characterisation by concentration-compactness techniques of the
nature of the possible blow-up in the critical case leading to the
global existence for this critical value. Finally, the last section is
devoted to the study of the free energy functional in self-similar
variables and the proof of the existence of self-similar solutions
in the sub-critical case.



\section{Existence criterion}\label{sec:existence}
 \setcounter{equation}{0}
As in \cite{S1,S2}, we consider the regularised problem
\begin{equation}\label{eq:spreg}
\left\lbrace
\begin{array}{rll}
\displaystyle \frac{\partial u_\eps}{\partial t}(t,x)&={\rm div}
\left[\nabla \left( f_\eps\!\circ\!u_\eps \right)(t,x)-u_\eps(t,x)\nabla
\phi_\eps(t,x)\right]\qquad & t>0\,,\;x\in\RR^d\;,\vspace{.3cm}\\
\displaystyle \phi_\eps(t,x)&=\EE \ast u_\eps(t,x)\;,\qquad
&t>0\,,\;x\in\RR^d\;,\vspace{.3cm}\\
u_\eps(0,x)&=u_0^\eps\geq 0\qquad &x\in\RR^d\,,
\end{array}\right.
\end{equation}
where $f_\eps:[0,\infty)\longrightarrow \RR$ is given by
$f_\eps(u):=(u+\eps)^m-\eps^m$. Here, $u_0^\eps$ is the convolution
of $u_0$ with a sequence of mollifiers and
$\|u_0^\varepsilon\|_1=\|u_0\|_1=M$ in particular. This regularised
problem has global in time smooth solutions. This approximation has been proved to be convergent. More
precisely, the result in \cite[Section 4]{S2} asserts that if we
assume that
\begin{equation}\label{ass:linfty}
\sup_{0<t<T} \|u_\eps(t)\|_\infty\leq \kappa
\end{equation}
where $\kappa$ is independent of $\eps>0$, then there exists a
sub-sequence $\eps_n\to 0$, such that
\begin{eqnarray}
 u_{\eps_n} \to u \qquad & \mbox{strongly} \qquad & \mbox{in }
 C([0,T],\LL^p_{{\rm loc}}(\RR^d)) \mbox{ and a.e. in } (0,T)\times\RR^d
,\label{conv1}\\
\nabla u_{\eps_n}^m \wto \nabla u^m \qquad & \mbox{weakly-*}
\qquad & \mbox{in } \LL^\infty(0,T;\LL^2(\RR^d)),\label{conv2}\\
\phi_{\eps_n}(t) \to \phi(t) \qquad & \mbox{strongly} \qquad &
\mbox{in } \LL^r_{{\rm loc}}(\RR^d) \mbox{ a.e. in } (0,T),\label{conv3}\\
\nabla \phi_{\eps_n}(t) \to \nabla \phi(t) \qquad & \mbox{strongly}
\qquad & \mbox{in } \LL^r_{{\rm loc}}(\RR^d) \mbox{ a.e. in }
(0,T),\label{conv4}
\end{eqnarray}
for any $p\in(1,\infty)$ and $r\in(1,\infty]$, and $u$ is a weak
solution to \eqref{eq:spf} on $[0,T)$ with $\phi=\EE\ast u$. Moreover, the free
energy \cite[Proposition 6.1]{S1} satisfies $\FF[u(t)]\leq
\FF[u_0]$ for a.e. $t\geq 0$. However, a detailed analysis of the
proof in \cite[Proposition 6.1]{S1} shows that the weak solution
is in fact a free energy solution.

\begin{proposition}[Existence of free energy
solutions]\label{prop:existencemaxentsol}
Under assumption \eqref{Eqn:Assumptions} on the initial data and
\eqref{ass:linfty} on the approximation sequence, there exists a
free energy solution to \eqref{eq:spf} in $[0,T)$.
\end{proposition}
\proof The only remaining points not covered by the results in
\cite{S1,S2} are the lower semi-continuity of the free energy
dissipation and the fact that $u^{(2m-1)/2}$ belongs to
$L^2(0,t;H^1(\RR^d))$ for $t\in [0,T)$. The latter will actually
be shown in the proof of Lemma~\ref{lem:extension}, see
\eqref{alet} below. Concerning the former, a
careful reading of the proof of \cite[Proposition 6.1]{S1} gives
that
$$
\FF_{\eps,l}[u_\eps(t)] + \frac34
\int_0^t\int_{\RR^d}[u_\eps(s,x)+\eps]\left|\nabla\left(\frac{m}{m-1}
[u_\eps(s,x)+\eps]^{m-1} -\phi_\eps(s,x) \right)\right|^2 \psi_l(x)
\dd x \dd t \leq \FF_{\eps,l}[u_0]
$$
for a.e. $t\in(0,T)$ where $\psi_l$ is a standard cut-off function in $\RR^d$ for any $l\in\NN$ and
$$
\FF_{\eps,l}[u_\eps(t)] = \int_{\RR^d} \frac{[u_\eps(t,x)+\eps]^m}{m-1}
\psi_l(x) \dd x - \frac12 \iint_{\RR^d \times \RR^d}
\EE(x-y)\,u_\eps(t,x)\,u_\eps(t,y)\dd x \dd y\;.
$$
In this regularised setting, we can write that
$$
(u_\eps+\eps)\left|\nabla\left[\frac{m}{m-1}
(u_\eps+\eps)^{m-1}-\phi_\eps\right]\right|^2 = \left| \frac{2m}{2m-1}\, \nabla
(u_\eps+\eps)^{(2m-1)/2} - (u_\eps+\eps)^{1/2}\nabla \phi_\eps
\right|^2\; .
$$
As proved in \cite{S1}, we have
$\FF_{\eps_n,l}[u_{\eps_n}(t)]\to\FF[u(t)]$ as $\eps_n\to 0$ and
$l\to\infty$. In addition, it is straightforward from the convergence
properties \eqref{conv1}-\eqref{conv4} above to pass to the limit as
$\eps_n\to 0$ in the free energy dissipation functional with the help
of a lower semi-continuity argument. We leave the details to the
reader, see {\it e.g.} \cite{otto-porous} or \cite[Lemma
10]{CJMTU}. Hence, passing to the limit as $l\to\infty$, then $u$ is a
free energy solution as it satisfies the free energy inequality
\eqref{freeeneinq}.
\finproof

\begin{remark}
The free energy inequality \eqref{freeeneinq} can be obtained with
constant $3/4$ multiplying the entropy dissipation directly from
{\rm \cite[Proposition 6.1]{S1}} and the procedure above. This is
a technical issue that can be improved to constant $1$ by redoing
the proof in {\rm \cite[Proposition~6.1]{S1}} treating more
carefully the free energy dissipation term. In fact, the proof in
{\rm \cite[Proposition~6.1]{S1}} shows that you can choose the
constant as close to $1$ as you want.
\end{remark}

We are now ready to characterise the maximal time of existence
by showing the local in time boundedness of the $\LL^m$-norm
independently of the approximation parameter $\eps>0$ and how this
estimate implies the local in time $\LL^\infty$-estimate
\eqref{ass:linfty}.


\begin{lemma}[From uniform integrability to $\LL^\infty$-bounds]
\label{lem:extension}
For any $\eta>0$ there exists $\tau_\eta>0$ depending only on $d$,
$M$, and $\eta$ such that, if
\begin{equation*}
\sup_{\eps\in (0,1)} \|u_\eps(t^*)\|_m \le \eta
\end{equation*}
for some $t^*\in [0,\infty)$, then

\begin{itemize}
\item[(i)] the family $(u_\eps)_\eps$ is bounded in
$\LL^\infty(t^*,t^*+\tau_\eta;\LL^m(\RR^d))$.

\item[(ii)] Moreover, if $(u_\eps(t^*))_\eps$ is also bounded in
$\LL^p(\RR^d)$ for some $p\in (m,\infty]$, then $(u_\eps)_\eps$ is
bounded in  $\LL^\infty(t^*,t^*+\tau_\eta;\LL^p(\RR^d))$.
\end{itemize}
\end{lemma}
\proof To prove this result we need to refine the argument already
used in the two-dimensional situation $d=2$ with linear diffusion
$m=1$ in \cite{BDP,BCM}. We follow a procedure analogous to the ones in
\cite{Kowalczyk04,Carrillo-Calvez,S1,CPZ}.

\medskip

\noindent{\bf Step~1 - $\LL^m$-estimates:} By \eqref{eq:spreg} we have
\begin{align*}
\frac{\dd}{\dd t} \|u_\eps\|_m^m =\,& -m(m-1) \int_{\RR^d}
u_\eps^{m-2}\,\nabla u_\eps\cdot \left( m (u_\eps +
\eps)^{m-1}\,\nabla u_\eps - u_\eps\,\nabla \phi_\eps \right) \dd x \\
\le\, & - \frac{4m^2(m-1)}{(2m-1)^2} \left\| \nabla u_\eps^{(2m-1)/2}
\right\|_2^2 - (m-1) \int_{\RR^d} u_\eps^m\,\Delta \phi_\eps \dd x \\
=\,& - \frac{4m^2(m-1)}{(2m-1)^2} \left\| \nabla u_\eps^{(2m-1)/2}
\right\|_2^2 + (m-1)\,\|u_\eps\|_{m+1}^{m+1}\; .
\end{align*}
As
$$
1 < \frac{2m}{2m-1} < \frac{2(m+1)}{2m-1} < \frac{2d}{d-2}\;,
$$
we have the following Gagliardo-Nirenberg-Sobolev inequality: there exists a positive constant $C$ such that
$$
\| w\|_{2(m+1)/(2m-1)} \le C\, \|\nabla
w\|_2^{[(2m-1)d]/[(m+1)(2m+d-2)]}
\,\|w\|_{2m/(2m-1)}^{2m^2/[(m+1)(2m+d-2)]}
$$
which we apply with $w=u_\eps^{(2m-1)/2}$ to obtain
\begin{equation*}
\|u_\eps\|_{m+1}^{(2m-1)/2} \le\, C\, \left\| \nabla u_\eps^{(2m-1)/2}
\right\|_2^{[(2m-1)d]/[(m+1)(2m+d-2)]}\|u_\eps\|_m^{m^2(2m-1)/[(m+1)(2m+d-2)]}\;.
\end{equation*}
It leads to
\begin{align*}
\|u_\eps\|_{m+1}^{m+1} \le\, & C\, \left\| \nabla u_\eps^{(2m-1)/2}
\right\|_2^{2d/(2m+d-2)}\,
\|u_\eps\|_m^{2m^2/(2m+d-2)} \\
\le \, & \frac{2m^2}{(2m-1)^2}\, \left\| \nabla u_\eps^{(2m-1)/2}
\right\|_2^2 + C\, \|u_\eps\|_m^{m^2/(m-1)}\;.
\end{align*}
We thus end up with
\begin{equation}
\label{vittel}
\frac{\dd}{\dd t} \|u_\eps\|_m^m + \frac{2m^2(m-1)}{(2m-1)^2}\,
\left\| \nabla u_\eps^{(2m-1)/2} \right\|_2^2 \le (m-1)\, C\,
\|u_\eps\|_m^{m^2/(m-1)}\;.
\end{equation}
In particular, for any $t_2\ge t_1\ge 0$
\begin{equation}
\label{luchon}
\|u_\eps(t_2)\|_m^m \le \left[ \|u_\eps(t_1)\|_m^{-m/(m-1)} -
C\,(t_2-t_1) \right]^{-(m-1)}
\end{equation}
Taking $t_1=t^*$ we deduce from \eqref{luchon} that
\begin{equation*}
\|u_\eps(t)\|_m^m \le \left( \eta^{-m/(m-1)} -
C\,(t-t^*) \right)^{-(m-1)}\quad \mbox{for $t\in [t^*,t^*+2\tau_\eta)$}
\end{equation*}
 with $\tau_\eta=1/\left(
2C\eta^{m/(m-1)} \right)$. Consequently, $\|u_\eps(t)\|_m^m \le
(C\tau_\eta)^{-(m-1)}$ for $t\in [t^*,t^*+\tau_\eta]$ and the proof
of the first assertion of Lemma~\ref{lem:extension} is complete. In
addition, coming back to \eqref{vittel}, we further deduce that
\begin{equation}
\label{alet}
\int_{t^*}^{t^*+\tau_\eta} \left\| \nabla
u_\eps^{(2m-1)/2} \right\|_2^2 \le C(t^*,\eta)\;.
\end{equation}

\noindent{\bf Step~2 - $\LL^p$-estimates, $p\in (m,\infty)$:} For
$t\in [t^*,t^*+\tau_\eta]$, $K\ge 1$, and $p>m$, we infer from
\eqref{eq:spreg} that
\begin{align*}
  \frac{\dd }{\dd t} \left\| (u_\eps(t)-K)_+ \right\|_p^p \leq \,&
-m\,p\,(p-1)\,\int_{\RR^d}  (u_\eps(t)-K)_+^{p-2}\,
(u_\eps+\eps)^{m-1}\,|\nabla u_\eps|^2 \dd x \\
&+ p\,(p-1)\, \int_{\RR^d}
\left[ (u_\eps(t)-K)_+^{p-1} + K\,(u_\eps(t)-K)_+^{p-2} \right]\,\nabla
u_\eps\cdot\nabla \phi_\eps \dd x \\
\le \,& -m\,p\,(p-1)\,\int_{\RR^d}  (u_\eps(t)-K)_+^{m+p-3}\,
|\nabla u_\eps|^2 \dd x \\
&- \int_{\RR^d} \left[
(p-1)\,(u_\eps(t)-K)_+^p +p\,K\,(u_\eps(t)-K)_+^{p-1}
\right]\,\Delta \phi_\eps \dd x \\
\le\, & -\frac{4\,m\,p\,(p-1)}{(m+p-1)^2} \left\|\nabla\left[
(u_\eps(t)-K)_+^{(m+p-1)/2}\right]\right\|_2^2 + {\rm (I)}
\end{align*}
with
$$
{\rm (I)} := p\,K^2\,\|(u_\eps(t)-K)_+\|_{p-1}^{p-1} +
(2p-1)\,K\,\|(u_\eps(t)-K)_+\|_p^{p}
+ (p-1)\, \|(u_\eps(t)-K)_+\|_{p+1}^{p+1} \;.
$$

We now use the following interpolation inequality
\begin{equation*}
\|w\|_{p+1}^{p+1} \le C(p)\, \left\| \nabla \left( w^{(m+p-1)/2}
\right) \right\|_2^2\, \|w\|_1^{2/d}
\end{equation*}
which is a consequence of the Gagliardo-Nirenberg-Sobolev and H\"older
inequalities (see, e.g., \cite[Lemma~3.2]{S2}) to obtain
\begin{eqnarray*}
{\rm (I)} & \le & (p-1)\,K^2\, \|(u_\eps(t)-K)_+\|_p^{p} + K^2\,\left| \{ x\,:\,
u_\eps(t,x)\ge K \} \right| + (2p-1)\,K\,\|(u_\eps(t)-K)_+\|_p^{p} \\
& + & C(p)\, \left\| \nabla \left[ (u_\eps(t)-K)_+^{(m+p-1)/2} \right]
\right\|_2^2\, \|(u_\eps(t)-K)_+\|_1^{2/d} \;.
\end{eqnarray*}
Noting that
$$
\|(u_\eps(t)-K)_+\|_1 \le \|u_\eps(t)\|_m\,
\left( \frac{\|u_\eps(t)\|_1}{K} \right)^{(m-1)/m}
$$
and recalling that $\|u_\eps(t)\|_1=M$ we conclude that
\begin{eqnarray*}
{\rm (I)} & \le & C(p)\, \frac{\|u_\eps(t)\|_m^{2/d}}{K^{2(m-1)/m\,d}}\,
\left\| \nabla \left[ (u_\eps(t)-K)_+^{(m+p-1)/2} \right] \right\|_2^2
\\
& + & K\, [2p-1+(p-1)\,K]\, \|(u_\eps(t)-K)_+\|_p^{p} + K\, M\;.
\end{eqnarray*}
By Step~1, we may choose $K=K_*$ large enough such that
$$
C\, \frac{\|u_\eps(t)\|_m^{2/d}}{K_*^{2(m-1)/m\,d}} \le
\frac{4\,m\,p\,(p-1)}{(m+p-1)^2}
$$
for all $t\in [t^*,t^*+\tau_\eta]$ and $\eps\in (0,1)$, hence
$$
{\rm (I)}\le \frac{4\,m\,p\,(p-1)}{(m+p-1)^2}\, \left\| \nabla \left[
(u_\eps(t)-K_*)_+^{(m+p-1)/2} \right] \right\|_2^2 + C(p,t^*,\eta)\,
\left[ 1 + \|(u_\eps(t)-K_*)_+\|_p^{p} \right]\;.
$$
Therefore
$$
\frac{\dd }{\dd t} \left\| (u_\eps(t)-K_*)_+ \right\|_p^p \le
C(p,t^*,\eta)\, \left[ 1 + \|(u_\eps(t)-K_*)_+\|_p^{p} \right]\;,
$$
so that
$$
\left\| (u_\eps(t)-K_*)_+ \right\|_p^p \le C(p,t^*,\eta) \;\;\mbox{
for }\;\; t\in [t^*,t^*+\tau_\eta] \;\;\mbox{ and }\;\; \eps\in
(0,1)\;.
$$
As
$$
\|u_\eps(t)\|_p^p \le C(p)\, \left( K_*^{p-m}\,\|u_\eps(t)\|_m^m +
\left\| (u_\eps(t)-K_*)_+ \right\|_p^p \right)\;,
$$
the previous inequality and Step~1 warrant that
$$
\left\| u_\eps(t) \right\|_p \le C(p,t^*,\eta) \;\;\mbox{ for
}\;\; t\in [t^*,t^*+\tau_\eta] \;\;\mbox{ and }\;\; \eps\in (0,1)\;.
$$

\medskip

\noindent{\bf Step~3 - $\LL^\infty$-estimates:} As a direct
consequence of Step~2 with $p=d+1$ and Morrey's embedding
theorem $(\nabla \phi_\eps)_\eps$ is bounded in
$\LL^\infty((t^*,t^*+\tau_\eta)\times\RR^d;\RR^d)$. This property in
turn implies that $(u_\eps)_\eps$ is bounded in
$\LL^\infty((t^*,t^*+\tau_\eta)\times\RR^d)$ and we refer to
\cite[Lemma 3.2]{Carrillo-Calvez} and \cite{Kowalczyk04} for a proof
(see also \cite[Section 5]{S2} and \cite[Theorem 1.2]{S1} for
alternative arguments). \finproof

\medskip

As a consequence of the previous lemma, we are able to construct a
free energy solution defined on a maximal existence time.
\begin{theorem}[Maximal free energy solution]\label{charac}
Under assumption \eqref{Eqn:Assumptions} on the initial condition there
are $T_\omega\in (0,\infty]$ and a free energy solution $u$ to
\eqref{eq:spf} on $[0,T_\omega)$ with the following alternative:
Either $T_\omega=\infty$ or $T_\omega<\infty$ and $\|u(t)\|_m\to\infty$
as $t\nearrow T_\omega$. Furthermore there exists a positive constant $C_0$ depending only on $d$ such that $u$ satisfies
\begin{equation}
  \label{eq:1}
  \|u(t_2)\|_m^m \le \left( \|u(t_1)\|_m^{-m/(m-1)} - C_0\,(t_2-t_1) \right)^{-(m-1)}
\end{equation}
for $t_1\in [0,T_\omega)$ and $t_2\in (t_1,T_\omega)$.
\end{theorem}
\proof We put $\xi_p(t)=\sup_{\eps\in (0,1)} \|u_\eps(t)\|_p \in
(0,\infty]$ for $t\ge 0$ and $p\in [m,\infty]$ and
$$
T_1 = \sup\left\{ T>0 \,: \, \xi_m\in \LL^\infty(0,T) \right\}\;.
$$
Clearly the definition of the sequence $(u_0^\eps)_\eps$ and~\eqref{Eqn:Assumptions} ensure that $\xi_p(0)$ is finite for all $p\in
[m,\infty]$. By Lemma~\ref{lem:extension} there exists $t_1>0$ such that
$\xi_p$ is bounded on $[0,t_1]$ for all $p\in [m,\infty]$. Then~\eqref{ass:linfty} is fulfilled for $T=t_1$ and there is a free
energy solution to~\eqref{eq:spf} on $[0,t_1)$ by
Proposition~\ref{prop:existencemaxentsol} and~\eqref{alet}. This
ensures in particular that $T_1\ge t_1>0$.

We next claim that
\begin{equation}
\label{quezac}
\xi_\infty \in \LL^\infty(0,T) \quad\mbox{for any $T\in[0,T_1)$}\;.
\end{equation}
Indeed, consider $T_1^\infty = \sup\{ T\in (0,T_1) \, : \,
\xi_\infty \in \LL^\infty(0,T) \}$ and assume for contradiction that
$T_1^\infty<T_1$. Then $\xi_m$ belongs to $\LL^\infty(0,T_1^\infty)$
and we put $\eta=\|\xi_m\|_{\LL^\infty(0,T_1^\infty)}$ and
$t^*=T_1^\infty-(\tau_\eta/2)$, $\tau_\eta$ being defined in
Lemma~\ref{lem:extension}. As $\xi_m(t^*)\le \eta$ and
$\xi_\infty(t^*)$ is finite we may apply
Lemma~\ref{lem:extension} to deduce that both $\xi_m$ and $\xi_\infty$
belong to $\LL^\infty(t^*,t^*+\tau_\eta)$, the latter property
contradicting the definition of $T_1^\infty$ as
$t^*+\tau_\eta=T_1^\infty+(\tau_\eta/2)$.

Now, thanks to \eqref{quezac}, \eqref{ass:linfty} is fulfilled for any
$T\in [0,T_1)$ and the existence of a free energy solution $u$ to
\eqref{eq:spf} on $[0,T_1)$ follows from
Proposition~\ref{prop:existencemaxentsol} and \eqref{alet}. Moreover, either
$T_1=\infty$ or $T_1<\infty$ and $\|u(t)\|_m\to\infty$ as
$t\nearrow T_1$, and the proof of Theorem~\ref{charac} is complete
with $T_\omega=T_1$. Or $T_1<\infty$ and
$$
\liminf_{t\to T_1} \|u(t)\|_m < \infty\;.
$$
In that case, there are $\eta>0$ and an increasing sequence of positive real
numbers $(s_j)_{j\ge 1}$ such that $s_j\to T_1$ as $j\to\infty$ and
$\|u(s_j)\|_m \le \eta$. Fix $j_0\ge 1$ such that $s_{j_0}\ge T_1 -
(\tau_\eta/2)$ with $\tau_\eta$ defined in Lemma~\ref{lem:extension}
and put $\tilde{u}_0=u(s_{j_0})$; According to
Definition~\ref{def:wfes} and \eqref{conv2} $\tilde{u}_0$ fulfils~\eqref{Eqn:Assumptions} and we may proceed as above to obtain a free
energy solution $\tilde{u}$ to \eqref{eq:spf} on $[0,T_2)$ for some
$T_2\ge\tau_\eta$. Setting $\bar{u}(t)=u(t)$ for $t\in [0,s_{j_0}]$
and $\bar{u}(t)=\tilde{u}(t-s_{j_0})$ for $t\in [s_{j_0},s_{j_0}+T_2)$
we first note that $\bar{u}$ is a free energy solution to
\eqref{eq:spf} on $[0,s_{j_0}+T_2)$ and a true extension of $u$ as
$s_{j_0}+T_2 \ge T_1 - (\tau_\eta/2)+\tau_\eta \ge T_1 +
(\tau_\eta/2)$. We then iterate this construction as long as the
alternative stated in Theorem~\ref{charac} is not fulfilled to complete
the proof.

Thanks to the regularity of weak solutions we may next proceed as
in the proof of~\eqref{luchon} to deduce~\eqref{eq:1}.
\finproof

\begin{corollary}[Lower bound on the blow-up time]
Let $u$ be a free energy solution
to~\eqref{eq:spf} on $[0,T_\omega)$ with an initial condition
$u_0$ satisfying~\eqref{Eqn:Assumptions}. If $T_\omega$ is finite,
then
\begin{equation*}
  \|u(t)\|_m \ge  \left[ C_0\,(T_\omega-t)\right]^{-(m-1)/m}\;,
\end{equation*}
where $C_0$ is defined in Theorem~{\rm\ref{charac}}.
\end{corollary}
\proof Let $t \in (0,T_\omega)$ and $t_2 \in (t,T_\omega)$. By ~\eqref{eq:1}, we have
\begin{equation*}
   \|u(t_2)\|_m^{-m/(m-1)} \ge \|u(t)\|_m^{-m/(m-1)} - C_0\,(t_2-t)\;.
\end{equation*}
Letting $t_2$ going to $T_\omega$ gives
\begin{equation*}
  0 \ge \|u(t)\|_m^{-m/(m-1)} - C_0\,(T_\omega-t)\,,
\end{equation*}
hence the expected result.
\finproof

\section{The free energy functional $\FF$}
 \setcounter{equation}{0}
As we have just seen in the existence proof, the existence time
of a free energy solution to \eqref{eq:spf} heavily depends on
the behaviour of its $\LL^m$-norm. As the free energy $\FF$ involves
the $\LL^m$-norm, the information given $\FF$ will be of paramount
importance. Let us then proceed to a deeper study of this functional.

\begin{lemma}[Scaling properties of the free
energy]\label{lem:scalingfreenrj} Given $h \in
\LL^1(\RR^d)\cap\LL^m(\RR^d)$, let us define
$h_\lambda(x):=\lambda^d h(\lambda\,x)$, then
 \begin{equation*}
   \FF[h_\lambda] = \lambda^{d-2} \FF[h]\quad \;\;\mbox{ for all $\lambda \in (0,\infty)$}\;.
 \end{equation*}
\end{lemma}
\proof We have
\begin{align*}
   \FF[h_\lambda] &= \frac{1}{m-1} \int_{\RR^d} \lambda^{2d-2}
h(\lambda\,x)^m \dd x - \frac{c_d}{2} \iint_{\RR^d\times\RR^d}
\lambda^{2d}\frac{h(\lambda\,x)\,h(\lambda\,y)}{|x-y|^{d-2}} \dd x \dd
y\\
&= \frac{\lambda^{d-2}}{m-1} \int_{\RR^d}  h(x)^m \dd x -
\frac{c_d\,\lambda^{d-2}}{2} \iint_{\RR^d\times\RR^d}
\frac{h(x)\,h(y)}{|x-y|^{d-2}} \dd x \dd y\\
&= \lambda^{d-2} \FF[h] \;,
 \end{align*}
giving the announced scaling property. \finproof

\medskip

We next establish a {\it variant to the Hardy-Littlewood-Sobolev
(VHLS) inequality}:
\begin{lemma}[VHLS inequality]\label{lem:vhls}
For $h \in \LL^1(\RR^d)\cap\LL^m(\RR^d)$ we put
$$
\WW(h) := \iint_{\RR^d\times\RR^d}\frac {h(x)\,h(y)}{|x-y|^{d-2}} \dd
x\dd y\,.
$$
Then
\begin{equation}\label{eq:adaptedHLSineq}
C_* := \,\sup_{h\neq 0}
\left\{\frac{\WW(h)}{\|h\|_1^{2/d}\,\|h\|_m^{m}}\,
, \, h\in\LL^1(\RR^d)\cap\LL^m(\RR^d)\right\} < \infty \ .
\end{equation}
\end{lemma}
First recall the Hardy-Littlewood-Sobolev (HLS) inequality, see \cite[Theorem~4.3]{liebloss}, which states that if
\begin{equation*}
 \frac 1p+\frac 1{q}+\frac\lambda{d}=2 \quad \mbox{and} \quad  0<\lambda<d\,,
\end{equation*}
then for all $f\in \LL^p(\R^d)$, $g\in \LL^q(\R^d)$, there exists
a sharp positive constant $C_{{\rm HLS}}>0$, given by~\cite{lieb83}, which only depends on $p$, $q$ and $\lambda$ such
that
\begin{equation}\label{eq:HLSineq}
\left|\iint_{\RR^d\times\RR^d}\frac {f(x)\,g(y)}{|x-y|^{\lambda}}
\dd x\dd y\right|\le C_{{\rm HLS}}\,\|f\|_{p}\,
\|g\|_{q}\;.
\end{equation}

\noindent {\sl Proof of Lemma~{\rm \ref{lem:vhls}}.\/}
Consider $h \in \LL^1(\RR^d)\cap\LL^m(\RR^d)$. Applying the HLS
inequality \eqref{eq:HLSineq} with $p=q=2d/(d+2)$ and $\lambda=d-2$, and then the H\"older inequality with $1<p=2d/(d+2)<m$, we obtain
$$
\left|\WW(h)\right| = \left|\iint_{\RR^d\times\RR^d}\frac
{h(x)\,h(y)}{|x-y|^{d-2}} \dd x\dd y\right|\le C_{{\rm HLS}}
\,\|h\|_p^2\leq C_{{\rm HLS}}
\,\|h\|_1^{2/d}\,\|h\|_m^{m}\ .
$$
Consequently, $C_*$ is finite and bounded from above by $C_{{\rm HLS}}$.
\finproof

\medskip

We next turn to the existence of maximisers for the VHLS inequality
which can be proved by similar arguments as for the classical HLS inequality in \cite[Theorem~2.5]{lieb83}.

\begin{lemma}[Extremals of the VHLS inequality]\label{lem:extremalfct}
There exists a non-negative, radially symmetric and non-increasing
function $P_* \in \LL^1(\RR^d)\cap\LL^m(\RR^d)$ such that
$\WW(P_*)=C_*$ with
$\|P_*\|_{1}=\|P_*\|_{m}=1$.
\end{lemma}
\proof Define
\begin{equation*}
  \Lambda(h):=\frac{\WW(h)}{\|h\|_1^{2/d}\,
\|h\|_m^m}
  \quad\mbox{ for $h \in \LL^1(\RR^d) \cap \LL^m(\RR^d)$}\;,
\end{equation*}
and consider a maximising sequence $(p_j)_j$ in $\LL^1(\RR^d) \cap
\LL^m(\RR^d)$, that is
\begin{equation}\label{eq:maxseq}
  \lim_{j \to \infty} \Lambda(p_j) = C_*\;.
\end{equation}\medskip

\noindent{\bf Step 1 -} We first prove that we may assume that
$p_j$ is a non-negative, radially symmetric, non-increasing
function such that
$\|p_j\|_1=\|p_j\|_m=1$ for any $j \ge
0$. Indeed, $\Lambda(p_j) \le \Lambda(|p_j|)$ so that $\left( |p_j|
\right)_j$ is also a maximising sequence. Next, let us introduce
$\tilde{p}_j(x):=\lambda_j |p_j(\mu_j x)|$ with
$\mu_j:=(\|p_j\|_1/\|p_j\|_m)^{m/[d(m-1)]}$ and
$\lambda_j:=\mu_j^d/\|p_j\|_1$. A
direct computation shows that
$\Lambda\left(\tilde{p}_j\right)=\Lambda(|p_j|)$ and
$\|\tilde{p}_j\|_1=\|\tilde{p}_j\|_m=1$.
Finally, denoting by $p_j^*$ the symmetric decreasing
rearrangement of $\tilde{p}_j$, we infer from the Riesz rearrangement properties~\cite[Lemma~2.1]{lieb83} that
\begin{equation*}
  \Lambda(p_j^*) = \WW(p_j^*)\ge  \WW(\tilde{p}_j) =
\Lambda(\tilde{p}_j) = \Lambda(|p_j|)\;.
\end{equation*}
Consequently, $\left( p_j^* \right)_j$ is also a maximising sequence
and the first step is proved.\medskip

\noindent{\bf Step 2 -}  Let us now prove that the supremum is
achieved. For $k \in \{1,m\}$, the monotonicity and the non-negativity
of $p_j$ imply that
\begin{equation*}
  1= \|p_j\|_k^k=d\,|B(0,1)|\int_0^\infty r^{d-1}p_j^k(r)\,\dd
r\ge d\,|B(0,1)|\,p_j^k(R)\int_0^R r^{d-1}\dd r \ge |B(0,1)|\,R^d\,p_j^k(R)\;.
\end{equation*}
So that
\begin{equation}
  \label{eq:15}
  0 \le p_j(R) \le b(R) := C_1\,\inf\{R^{-d/m};\,R^{-d}\}\quad \mbox{
for }\; R>0\;.
\end{equation}

Now, we use once more the monotonicity of the $p_j$'s and their
boundedness in $(R,\infty)$ for any $R>0$ to deduce from Helly's theorem that
there are a sub-sequence of $(p_j)_j$ (not relabelled) and a
non-negative and non-increasing function $P_*$ such that $(p_j)_j$
converges to $P_*$ point-wisely. In addition, as $1<2d/(d+2)<m$,
$x\mapsto b(|x|)$ belongs to $\LL^{2d/(d+2)}(\RR^d)$ while the
HLS inequality \eqref{eq:HLSineq} warrants
that
\begin{equation*}
  (x,y) \mapsto b(|x|)\,b(|y|)\,|x-y|^{-(d-2)} \in \LL^1 (\RR^d \times
\RR^d)\;.
\end{equation*}
Together with \eqref{eq:15} and the point-wise convergence of
$(p_j)_j$, this implies that
\begin{equation*}
  \lim_{j \to \infty} \WW(p_j)= \WW(P_*)
\end{equation*}
by the Lebesgue dominated convergence theorem. Consequently,
$\WW(P_*)=C_*$ and thus $P_* \neq 0$. In addition, the point-wise
convergence of $(p_j)_j$ and Fatou's lemma ensure $\|P_*\|_1 \le 1$
and $\|P_*\|_m \le 1$. Therefore $\Lambda(P_*) \ge C_*$ and using~\eqref{eq:maxseq} we conclude that $\Lambda(P_*) = C_*$. This in turn
implies that $\|P_*\|_1=\|P_*\|_m= 1$. \finproof

\medskip

We are now in a position to begin the study of the free energy
functional $\FF$. To this end, let us define the \textit{critical
mass} $M_c$ by 
\begin{equation}
  \label{eq:mc}
  M_c := \left[\frac{2}{(m-1)\,C_{*}\,c_d} \right]^{d/2}\;.
\end{equation}
Next, for $M>0$, we put
\begin{equation*}
 \mu_M:= \inf_{h \in \YY_{M}}\FF[h] \quad\mbox{where}\quad \YY_M:=\{ h
\in \LL^1(\RR^d) \cap \LL^m(\RR^d) \;:\; \|h\|_1=M\} \;,
\end{equation*}
and first identify the values of $\mu_M$ as a function of $M>0$.
\begin{proposition}[Infimum of the free energy]\label{prop:inffreenrj}
We have
\begin{equation}
\label{eq:27} \mu_M=\left\{
    \begin{array}{ll}
      0\quad &\mbox{if $M \in (0,M_c]$,}\vspace{.3cm}\\
-\infty \quad &\mbox{if $M >M_c$.}
    \end{array}
\right.
  \end{equation}
Moreover,
\begin{equation}\label{eq:lmbound}
\frac{C_{*}\,c_d}{2}\ \left( M_c^{2/d} - M^{2/d} \right)\
\|h\|_m^m \le \FF[h] \le \frac{C_{*}\,c_d}{2}\ \left(
M_c^{2/d} + M^{2/d} \right)\ \|h\|_m^m
\end{equation}
for $h\in\YY_M$. Furthermore, the infimum $\mu_M$ is not achieved if
$M<M_c$ while there exists one minimiser of $\FF$ in $\YY_{M_c}$.
\end{proposition}
\proof
Consider $h \in \LL^1(\RR^d) \cap \LL^m(\RR^d)$.
By the VHLS inequality \eqref{eq:adaptedHLSineq},
\begin{equation*}
\FF[h] \ge \left( \frac{1}{m-1} - \frac{C_{*}\,c_d}{2}
M^{2/d}\right)\|h\|^m_m\ge
\frac{C_{*}\,c_d}{2}\left( M_c^{2/d} -
M^{2/d}\right)\|h\|^m_m\;,
\end{equation*}
and
\begin{equation*}
\FF[h] \le \left( \frac{1}{m-1} + \frac{C_{*}\,c_d}{2}
M^{2/d}\right)\|h\|^m_m\le
\frac{C_{*}\,c_d}{2}\left( M_c^{2/d} +
M^{2/d}\right)\|h\|^m_m\;,
\end{equation*}
hence \eqref{eq:lmbound}.

\medskip

\noindent{\bf Case $M \le M_c\ $ -} By~\eqref{eq:lmbound}, $\FF$ is non-negative, so that $\mu_M\ge 0$. Choosing
\begin{equation*}
  h_*(t,x)= \frac{M}{(2\,\pi\,t)^{d/2}} \,
  \ee^{-|x|^2/(4t)}\;,
\end{equation*}
then
\begin{equation*}
\|h_*(t)\|_1=M
\quad\mbox{and}\quad\|h_*(t)\|^m_m=O\left(t^{-d(m-1)/2} \right)\;.
\end{equation*}
Therefore $h_*(t)$ belongs to $\YY_M$ for each $t>0$ and it follows
from \eqref{eq:lmbound} that $\FF[h_*(t)]\rightarrow 0$ as
$t\to\infty$. The infimum $\mu_M$ of $\FF$ on $\YY_M$ is thus
non-positive, hence $\mu_M=0$.

Finally, in the case $M<M_c$, $\mu_M=0$ and
\eqref{eq:lmbound} imply that the infimum of $\FF$ in $\YY_M$
is not achieved. If $M=M_c$ and $p\in \LL^1(\RR^d)\cap\LL^m(\RR^d)$ satisfies
$\WW(p)=C_{*}\,\|p\|_m^{m}\,\|p\|_1^{2/d}$ (such a
function exists by Lemma~\ref{lem:extremalfct}), then
\begin{equation*}
  \tilde{p}(x):=M_c^{-d/(d-2)} p\left(x M_c^{-m/(d-2)}\right)
\end{equation*}
belongs
to $\YY_{M_c}$ with $\|\tilde{p}\|_m=1$ and
$\WW(\tilde{p})=C_*\,M_c^{2/d}$. Therefore, $\FF[\tilde{p}]=0$ and we
have thus proved that suitably rescaled extremals of the VHLS
inequality \eqref{eq:adaptedHLSineq} are minimisers for $\FF$ in
$\YY_{M_c}$.

\medskip

\noindent{\bf Case $M > M_c\ $ -} This part of the proof is based
on arguments in \cite{Weinstein83}. Fix $\theta \in
\left((M_c/M)^{2/d},1 \right)$. By the VHLS inequality
\eqref{eq:adaptedHLSineq}, there exists a non-zero function $h^* \in
\LL^1(\RR^d) \cap \LL^m(\RR^d)$, such that
\begin{equation}
\label{eq:thetec}
  \theta \, C_* \le
\frac{|\WW(h^*)|}{\|h^*\|_m^{m}\|h^*\|_1^{2/d}} \le
C_*\;.
\end{equation}
Since $|\WW(h^*)|\le \WW(|h^*|)$ we may assume without loss of
generality that $h^*$ is non-negative. Let $\lambda>0$ and
consider the function $h_\lambda(x):=\lambda^d
h^*\left(\lambda\,\|h^*\|_1^{1/d}\,M^{-1/d} x\right)$. Then,
$h_{\lambda} \in \YY_{M}$ and it follows from the definition of
$M_c$ and \eqref{eq:thetec} that
\begin{align*}
  \FF[h_{\lambda}]&= \lambda^{d-2}
\left[\frac{M}{(m-1)\,\|h^*\|_1}\,\|h^*\|_m^m -
\frac{c_d}{2}\,\left( \frac{M}{\|h^*\|_1}
\right)^{(d+2)/d}\,\WW(h^*)\right]\\
& \le \lambda^{d-2}
\left[\frac{M}{(m-1)\,\|h^*\|_1}\,\|h^*\|_m^m -
\frac{c_d}{2}\,\left( \frac{M}{\|h^*\|_1}
\right)^{(d+2)/d}\,\theta\,C_*\,\|h^*\|_m^m\,
\|h^*\|_1^{2/d}\right]\\
& = \lambda^{d-2}\,\left( \frac{M}{\|h^*\|_1}
\right)^{(d+2)/d}\,\frac{\|h^*\|_m^m}{m-1}\,\left[ \left(
\frac{M_c}{M}\right)^{2/d}-\theta \right]\;.
\end{align*}
Owing to the choice of $\theta$ we may let $\lambda$ go to
infinity to obtain that $\mu_M=-\infty$, thus completing the proof.
\finproof

Let us now describe the set of minimisers of $\FF$ in $\YY_{M_c}$.

\begin{proposition}[Identification of the minimisers]\label{prop:caracterisationminimiser}
Let $\zeta$ be the unique positive radial classical solution to
$$
\Delta \zeta + \frac{m-1}{m}\, \zeta^{1/(m-1)}=0 \;\;\mbox{ in }\;\;
B(0,1) \;\;\mbox{ with }\;\; \zeta=0 \;\;\mbox{ on }\;\;\partial
B(0,1)\,.
$$
If $V$ is a minimiser of $\FF$ in $\YY_{M_c}$ there are $R>0$ and $z
\in \RR^d$ such that
  \begin{equation*}
   V(x)= \left\{
      \begin{array}{ll}
        \displaystyle \frac{1}{R^d} \left[\zeta\left(\frac{x-z}{R}
\right) \right]^{d/(d-2)} \quad & \mbox{if $x \in B(z,R)$,}\vspace{.3cm}\\
\displaystyle 0 \quad & \mbox{if $x \in \RR^d\setminus B(z,R)$.}
      \end{array}
 \right.
  \end{equation*}
\end{proposition}
\proof We have already shown in Proposition~\ref{prop:inffreenrj}
that the function $\FF$ has at least a minimiser in $\YY_{M_c}$. Let $V$ be a minimiser of $\FF$ in $\YY_{M_c}$, and define
$\tilde V(x):=\|V\|_m^{-m/(m-1)}V\left(x\,
\|V\|_m^{-m/(d(m-1))}\right)$ for $x \in \RR^d$. We have
$\|\tilde V\|_1=M_c$, $\|\tilde V\|_m=1$ and
$\FF[\tilde V]=0$, so
that $\tilde V$ is also a minimiser of $\FF$ in $\YY_{M_c}$. We
next denote by $W$ the symmetric rearrangement of $\tilde V$. Then
$\|W\|_1=M_c$, $\|W\|_m=\|\tilde V\|_m=1$ and
$\WW(W) \ge |\WW(\tilde V)|$ by the Riesz rearrangement properties~\cite[Lemma~2.1]{lieb83}.
Therefore, $\FF[W]\le \FF[\tilde
V]=0$ and thus $\FF[W]=0$ since $W \in \YY_{M_c}$. This in turn
implies that $\WW(W)=|\WW(\tilde V)|$. Again by \cite[Lemma~2.1]{lieb83}
there is $y \in \RR^d$ such that $\tilde V(x)=W(x-y)$ for $y \in
\RR^d$.

We next derive the Euler-Lagrange equation solved by $W$ and first
point out that a difficulty arises from the non-differentiability
of the $\LL^1$-norm. Nevertheless, we introduce $\Sigma_0:=\{x \in
\RR^d\,:\, W(x)=0\}$, $\Sigma_+:=\{x \in \RR^d\,:\, W(x)>0\}$ and
consider $\varphi \in \CC^\infty_0(\RR^d)$ and $\eps>0$. The
perturbation
$M_c\,\|W+\eps\,\varphi\|_1^{-1}(W+\eps\,\varphi)$ belongs
to $\YY_{M_c}$ and is such that
\begin{equation*}
  \FF\left[
\frac{M_c}{\|W+\eps\,\varphi\|_1}(W+\eps\,\varphi)\right] \ge
\FF[W] \ge 0\;.
\end{equation*}
After a few computations that we omit here we
may let $\eps \to 0$, and conclude that
\begin{align*}
  &2\iint_{\RR^d \times \RR^d}
\frac{W(x)\,\varphi(y)}{|x-y|^{d-2}}\dy\dx \\
&\hspace{1cm}\le C_*\,M_c^{2/d}\,m\int_{\RR^d}
W^{m-1}(x)\varphi(x)\dx +\frac{2}{d} \,C_*\,M_c^{(2-d)/d}
\left( \int_{\Sigma_+} \varphi(x)\dx+ \int_{\Sigma_0}
|\varphi(x)|\dx\right)\;.
\end{align*}
Using the definition of $M_c$ and $\EE$, the above formula also reads
\begin{equation}
\label{spirou}
\int_{\RR^d} \left( \frac{m}{m-1}\, W^{m-1} - \EE \ast W +
\frac{2-m}{m-1}\, \frac{1}{M_c} \right)\,\varphi\dx \ge
\frac{2-m}{m-1}\, \frac{1}{M_c}\, \int_{\Sigma_0}
(\varphi-|\varphi|)\dx
\end{equation}
for all $\varphi\in \CC^\infty_0(\RR^d)$. On the one hand, the
right-hand side of \eqref{spirou} vanishes for any non-negative
$\varphi\in \CC^\infty_0(\RR^d)$, so that
$$
\frac{m}{m-1}\, W^{m-1} - \EE \ast W + \frac{2-m}{m-1}\,
\frac{1}{M_c} \ge 0 \;\;\mbox{ a.e. in }\;\;\RR^d\,.
$$
Therefore, for almost every $x\in\Sigma_0$, we have $0 \ge \EE \ast
W(x) - (2-m)/[(m-1)M_c]$ so that
\begin{equation}
\label{fantasio}
\frac{m}{m-1}\, W^{m-1}(x) = 0 = \left( \EE \ast W(x) -
\frac{2-m}{m-1}\, \frac{1}{M_c} \right)_+ \;\;\mbox{ for almost every
}\;\;  x\in\Sigma_0\,.
\end{equation}
On the other hand, if $\psi\in \CC^\infty_0(\RR^d)$, a standard
approximation argument allows us to take
$\varphi=\mathbf{1}_{\Sigma_+}\,\psi$ in \eqref{spirou} and deduce
that
$$
\int_{\Sigma_+} \left( \frac{m}{m-1}\, W^{m-1} - \EE \ast W +
\frac{2-m}{m-1}\, \frac{1}{M_c} \right)\,\psi\dx \ge 0\,.
$$
This inequality being also valid for $-\psi$, we conclude that the
left-hand side of the above inequality vanishes for all $\psi\in
\CC^\infty_0(\RR^d)$, whence
\begin{equation}
\label{spip}
\frac{m}{m-1}\, W^{m-1} = \EE \ast W - \frac{2-m}{m-1}\,
\frac{1}{M_c} \;\;\mbox{ a.e. in }\;\;\Sigma_+\,.
\end{equation}
Combining \eqref{fantasio} and \eqref{spip} gives
\begin{equation*}
\label{gaston}
\frac{m}{m-1}\, W^{m-1} = \left( \EE \ast W - \frac{2-m}{m-1}\,
\frac{1}{M_c} \right)_+ \;\;\mbox{ a.e. in }\;\;\RR^d\,.
\end{equation*}
Now, since $W$ is radially symmetric and non-increasing there exists
$\rho \in (0,\infty]$ such that
\begin{equation*}
  \Sigma_+ \subset B(0,\rho)\;\mbox{ and }\;\Sigma_0 \subset
\RR^d\setminus B(0,\rho)\;,
\end{equation*}
and we infer from \eqref{spip} that
\begin{equation}
  \label{eq:22}
 \frac{m}{m-1} W^{m-1} = \EE \ast W - \frac{2-m}{m-1}\frac{1}{M_c}
\quad \mbox{for a.e. $x \in B(0,\rho)$}\;.
\end{equation}
Since $W \in \LL^r(\RR^d)$ for each $r \in (1,m]$ it follows from
the HLS inequality~\eqref{eq:HLSineq} that $\EE \ast W \in
\LL^r(\RR^d)$ for each $r\in \left(d/(d-2),m/(m-1)^2 \right]$,
see \cite[Theorem~10.2]{liebloss}. In particular, $\EE \ast W$ and
$W^{m-1}$ both belong to $\LL^{m/(m-1)}(\RR^d)$. This property and
\eqref{eq:22} then exclude that $\rho=\infty$ as $M_c>0$. Therefore
$\rho<\infty$ and
 \begin{equation*}
   \frac{m}{m-1}\, W^{m-1}(x)=
\left\{
   \begin{array}{ll}
     \displaystyle \EE \ast W(x) -
\frac{2-m}{m-1}\frac{1}{M_c}&\;\mbox{if $|x| <\rho$}\;,\vspace{.3cm}\\
     0&\;\mbox{if $|x|  > \rho$}\;.
   \end{array}
\right.
 \end{equation*}
Since $\EE \ast W \in \LL^{m/(m-1)^2}(\RR^d)$, the above
inequality allows us to conclude that $W \in
\LL^{m/(m-1)}(\RR^d)$. We now improve the regularity of $W$ by classical elliptic estimates. Introduce $\theta:=W^{m-1}$
and note that
\begin{equation*}
  \frac{m}{m-1}\, \theta(x)= \int_{\RR^d}\EE(x-y)\,W(y)\dy +
\frac{m-1}{m-2}\frac{1}{M_c}
\end{equation*}
for $x \in B(0,\rho)$ and $W \in \LL^{m/(m-1)}(\RR^d)$. By~\cite[Theorem~9.9]{gilbargtrudinger}, we have $\theta \in
W^{2,m/(m-1)}(B(0,\rho))$. A bootstrap argument then ensures
that $\theta$ and $W$ both belong to $W^{2,r}(B(0,\rho))$ for every
$r\in (1,\infty)$. It then follows from
\cite[Lemma~4.2]{gilbargtrudinger} that $\theta \in \CC^2(B(0,\rho))$
with $-\Delta \theta = (m-1)\theta^{m/(m-1)}/m$ in
$B(0,\rho)$ while \cite[Lemma~4.1]{gilbargtrudinger} warrants that
$\theta \in \CC^1(\RR^d)$. Then $\theta(x)=0$ if $|x|=\rho$ and
$\theta$ is thus a classical solution to $-\Delta \theta =
(m-1)\theta^{m/(m-1)}/m$ in $B(0,\rho)$ with $\theta = 0$
on $\partial B(0,\rho)$. By~\cite[Lemma~2.3]{gidasninirenberg}, there
is a unique positive solution to this problem. In fact, a simple scaling
argument shows that
\begin{equation*}
   \theta(x)=\rho^{2(m-1)/(m-2)}\zeta\left(\frac{x}{\rho}
\right)\quad\mbox{ for $x$ in $B(0,\rho)$}
 \end{equation*}
 and then
 \begin{equation*}
   W(x)=\frac{1}{\rho^d}\left[ \zeta\left(\frac{x}{\rho}
\right)\right]^{d/(d-2)}\quad\mbox{ for $x$ in $B(0,\rho)$}\;.
 \end{equation*}
Coming back to $V$, we have
 \begin{equation*}
   V(x)=
 \left\{
   \begin{array}{ll}
     \displaystyle \lambda^d W(\lambda\,x-y)=0 &\mbox{ if }
\displaystyle x \in \RR^d\setminus
B\left(\frac{y}{\lambda},\frac{\rho}{\lambda} \right)\;,\vspace{.3cm}\\
\displaystyle  \left(\frac{\lambda}{\rho} \right)^d
\left[\zeta\left(\left(x-\frac{y}{\lambda}\right)\left(
\frac{\rho}{\lambda}\right)^{-1}\right) \right]^{d/(d-2)}&\mbox{ if }
\displaystyle x \in B\left(\frac{y}{\lambda},\frac{\rho}{\lambda}
\right)\;,
 \end{array}
 \right.
 \end{equation*}
 which is the desired result with $R=\rho/\lambda$ and $z=y/\lambda$.
\finproof
\begin{remark}
  As a consequence of the identification of the minimisers given in Proposition~{\rm \ref{prop:caracterisationminimiser}}, $C_*<C_{\rm HLS}$. Otherwise any minimiser $V$ of $\FF$ is $\YY_{M_c}$ would also be an extremum for the HLS inequality~\eqref{eq:HLSineq} and thus be equal to
  \begin{equation*}
    V(x)=\frac{a}{(1+|x|^2)^{(d+2)/2}}\;,
  \end{equation*}
 for some $a>0$, see {\rm \cite[Theorem~3.1]{lieb83}}. This contradicts Proposition~{\rm \ref{prop:caracterisationminimiser}}.
\end{remark}
\begin{lemma}[Unboundedness of $\FF$]\label{lem:supfreenrj}
For each $M>0$ we have
\begin{equation}
\label{eq:supfreenrj}
\sup_{h\in \YY_M} \FF[h] = +\infty\,.
\end{equation}
\end{lemma}
If $M\in (0,M_c)$ the claim \eqref{eq:supfreenrj} is actually a
straightforward consequence of~\eqref{eq:lmbound}.\smallskip 

\proof Let $M>0$ and assume for contradiction that
$$
A:=\sup_{h\in \YY_M} \FF[h]<\infty\,.
$$
Consider $h\in\LL^1(\RR^d)\cap\LL^m(\RR^d)$ and define $h_\lambda(x):= M\lambda^d h(\lambda
x)/\|h\|_1$ for $x\in\RR^d$ and $\lambda>0$. Then $\|h_\lambda\|_1=M$ so that $h_\lambda\in\YY_M$
with
$$
\|h_\lambda\|_m^m=\lambda^{d-2}\, \left( \frac{M}{\|h\|_1} \right)^m\, \|h\|_m^m \;\;\mbox{ and
}\;\; \WW(h_\lambda)=\lambda^{d-2}\, \left( \frac{M}{\|h\|_1} \right)^2\, \WW(h)\,.
$$
Since $h_\lambda\in\YY_M$ we have $\FF[h_\lambda]\le A$, hence
\begin{equation*}
\|h_\lambda\|_m^m  \le (m-1)\left( A + \frac{c_d}{2}\, \WW(h_\lambda)\right)\\
\end{equation*}
from which we deduce
\begin{align*}
\|h\|_m^m & \le (m-1)\left(A\, \lambda^{2-d}\, \left( \frac{\|h\|_1}{M} \right)^m +
\frac{c_d}{2}\, \left( \frac{M}{\|h\|_1} \right)^{2/d}\, \WW(h)\right)\\
& \le (m-1)A\, \lambda^{2-d}\, \left( \frac{\|h\|_1}{M} \right)^m + \frac{1}{C_*
M_c^{2/d}}\, \left( \frac{M}{\|h\|_1} \right)^{2/d}\, \WW(h)\,.
\end{align*}
This inequality being valid for all $\lambda>0$ we let $\lambda\to\infty$ and use the HLS
inequality \eqref{eq:HLSineq} to obtain
$$
\|h\|_m^m \le \frac{1}{C_* M_c^{2/d}}\, \left( \frac{M}{\|h\|_1} \right)^{2/d}\, \WW(h) \le \left(
\frac{M}{M_c} \right)^{2/d}\, \frac{C_{\rm HLS}}{C_*}\, \frac{\|h\|_{2d/(d+2)}^2}{\|h\|_1^{2/d}}\,.
$$
Consequently,
\begin{equation}\label{eq:reverseholder}
\|h\|_m^m\, \|h\|_1^{2/d} \le \left( \frac{M}{M_c} \right)^{2/d}\, \frac{C_{\rm HLS}}{C_*}\,
\|h\|_{2d/(d+2)}^2
\end{equation}
for all $h\in\LL^1(\RR^d)\cap\LL^m(\RR^d)$.

Now, as $2d/(d+2)<m$, we may choose $\gamma\in ((d+2)/d,d/m)$ and put
$b_\delta(x):=(\delta+|x|)^{-\gamma}\, \un_{B(0,1)}(x)$ for $x\in\RR^d$ and $\delta\in [0,1]$.
Clearly $b_\delta$ belongs to $\LL^1(\RR^d)\cap\LL^m(\RR^d)$ with $\|b_\delta\|_1\ge \|b_1\|_1>0$
and $\|b_\delta\|_{2d/(d+2)}\le \|b_0\|_{2d/(d+2)}<\infty$ for each $\delta\in (0,1]$. These
properties and \eqref{eq:reverseholder} readily imply that $(b_\delta)_{\delta\in (0,1]}$ is
bounded in $\LL^m(\RR^d)$ which is clearly not true according to the choice of $\gamma$. Therefore
$A$ cannot be finite and Lemma~\ref{lem:supfreenrj} is proved.
\finproof
\section{Critical threshold}
 \setcounter{equation}{0}
It turns out that the critical mass $M_c$ arising in the study of the
free energy functional and defined in \eqref{eq:mc} plays also an
important role in the dynamics of \eqref{eq:spf}. In the next sections
we will distinguish the three cases $M>M_c$ (\textit{super-critical case}),
$M<M_c$ (\textit{sub-critical case}), and $M=M_c$ (\textit{critical
case}), $M$ denoting the $L^1(\RR^d)$-norm of the initial condition
$u_0$.

\subsection{Finite time blow-up in the super-critical case}

We start with the case $M>M_c$ in which we use the standard
argument relying on the evolution of the second moment of solutions as originally done
in~\cite{Jager92} for the PKS system corresponding to $d=2$
and $m=1$.
\begin{lemma}[Virial identity]\label{lem:viriel}
Under assumption \eqref{Eqn:Assumptions}, let $u$ be a free energy
solution to \eqref{eq:spf} on $[0,T)$ with initial condition $u_0$
for some $T\in (0,\infty]$. Then
\begin{equation*}
   \frac{\dd}{\dd t} \int_{\RR^d} |x|^2\,u(t,x)\dd x =
2\,(d-2)\,\FF[u(t)]\,,\quad t\in [0,T)\;.
\end{equation*}
\end{lemma}
\proof Here, we show the formal computation leading to this
property, the passing to the limit from the approximated problem
\eqref{eq:spreg} can be done by adapting the arguments in
\cite[Lemma~6.2]{S1} and \cite[Lemma~2.1]{BDP} without any further
complication. By integration by parts in \eqref{eq:spf} and
symmetrising the second term, we obtain
\begin{align*}
\frac{\dd}{\dd t} \int_{\RR^d} |x|^2\,u(t,x)\dd x&=
2\,d\,\int_{\RR^d}u^m(t,x)\dd x + 2 \iint_{\RR^d\times\RR^d}
[x\cdot\nabla \EE (x-y)]\,u(t,x)\,u(t,y) \dd y \dd x\\
&= 2\,d\,\int_{\RR^d}u^m(t,x)\dd x + \iint_{\RR^d\times\RR^d}
[(x-y)\cdot\nabla \EE (x-y)]\,u(t,x)\,u(t,y) \dd y \dd x \\
&=2\,(d-2) \,\FF[u(t)]\;,
\end{align*}
giving the desired identity. \finproof

Let us mention that a similar argument can be found in \cite[Lemma 6.2]{S1} and ~\cite{S2} in the present situation where the Poisson equation is
substituted by $-\Delta \phi= u - \phi$. The previous evolution for the second moment is
simpler in our case than the one in~\cite{S2} and resembles
that arising in the study of critical nonlinear Schr\"odinger
equations~\cite{Ca03}.

Let us also emphasise that this second moment evolution is more complicated than in the classical PKS system corresponding to $d=2$ and $m=1$ where the time derivative of the second moment is a constant. \medskip

An easy consequence of the previous lemma is the following
blow-up result.

\begin{proposition}[Blowing-up solutions]\label{lem:superfinite}
If $M>M_c$, then there are initial data $u_0$ satisfying~\eqref{Eqn:Assumptions} with $\|u_0\|_1=M$ and
negative free energy $\FF[u_0]$. Moreover, if $u_0$ is such
an initial condition and $u$ denotes a free energy
solution to \eqref{eq:spf} on $[0,T_\omega)$ with initial condition
$u_0$, then $T_\omega<\infty$ and the $\LL^m$-norm of $u$ blows up
in finite time.
\end{proposition}
\proof The proof is based on the idea of Weinstein
\cite{Weinstein83}. By the identification of the minimisers for the critical mass given in
Proposition~\ref{prop:caracterisationminimiser},
$\tilde{u}:=\zeta^{d/(d-2)}$ satisfies \eqref{Eqn:Assumptions} as well as
$\|\tilde{u}\|_1=M_c$ and $\FF[\tilde{u}]=0$. For $M>M_c$, the initial
condition $u_0=(M/M_c)\tilde{u}$ also satisfies \eqref{Eqn:Assumptions} with
$\|u_0\|_1=M$ and
\begin{align*}
  \FF[u_0]&=\frac{1}{m-1} \left(\frac{M}{M_c}
\right)^m\|\tilde{u}\|_m^m - \left(\frac{M}{M_c} \right)^2
\frac{c_d}{2} \WW(\tilde{u})\\
&=\frac{1}{m-1}\left[ \left(\frac{M}{M_c} \right)^m -
\left(\frac{M}{M_c} \right)^2\right]\|\tilde{u}\|_m^m\;,
\end{align*}
is negative as $M>M_c$ and $m<2$.

Consider next an initial condition $u_0$ satisfying
\eqref{Eqn:Assumptions} as well as $\|u_0\|_1>M_c$ and $\FF[u_0]<0$.
Denoting by $u$ a corresponding free energy solution to
\eqref{eq:spf} on $[0,T)$, we infer from the time
monotonicity of $\FF$ and Lemma~\ref{lem:viriel} that
\begin{equation*}
  \frac{\dd}{\dd t} \int_{\RR^d} |x|^2\,u(t,x)\dd x =
2\,(d-2)\,\FF[u(t)] \le  2\,(d-2)\,\FF[u_0] <0\;.
\end{equation*}
This implies that the second moment of $u(t)$ will become negative
after some time and contradicts the non-negativity of $u$. Therefore,
$T_\omega$ is finite and $\|u\|_m$ blows up in finite time.\finproof
\subsection{Global existence}
\begin{proposition}[Global existence in the subcritical
case]\label{prop:existenceglobalsol}
Under assumption \eqref{Eqn:Assumptions}, there exists a free energy
solution to \eqref{eq:spf} in $[0,\infty)$ with initial condition $u_0$.
\end{proposition}
\proof By Theorem~\ref{charac} there are $T_\omega$ and a free energy
solution to \eqref{eq:spf} in $[0,T_\omega)$ with initial condition $u_0$.
We then infer from \eqref{freeeneinq}, \eqref{evian}, and \eqref{eq:lmbound}
that $u(t)$ belongs to $\YY_M$ for all $t\in [0,T_\omega)$ and
$$
\frac{C_{*}\,c_d}{2}\ \left( M_c^{2/d} - M^{2/d} \right)\
\|u(t)\|_m^m \le \FF[u(t)] \le \FF[u_0]\;.
$$
As $M<M_c$, we deduce from the previous inequality that $u$ lies in
$\LL^\infty(0,\min{\{T,T_\omega\}}; \LL^m(\RR^d))$ for every $T>0$
which implies that $T_\omega=\infty$ by Theorem~\ref{charac}.
\finproof
\medskip

Let us now discuss the critical case.
\subsubsection{How would it blow-up?}
\begin{proposition}[Nature of the blow-up]\label{prop:how}
Let $u_0$ be an initial condition satisfying
\eqref{Eqn:Assumptions} with $\|u_0\|_1=M_c$ and consider a free
energy solution $u$ to \eqref{eq:spf} on $[0,T_\omega)$ with
initial condition $u_0$ and $T_\omega \in (0,\infty]$ and such
that $\|u(t)\|_m\to\infty$ as $t\nearrow T_\omega$. If $(t_k)_k$
is a sequence of positive real numbers such that $t_k\to T_\omega$
as $k\to\infty$, there are a sub-sequence $(t_{k_j})_j$ of
$(t_k)_k$ and a sequence $(x_j)_j$ in $\RR^d$ such that
\begin{equation*}
  \lim_{j \to \infty} \int_{\RR^d}\left| u(t_{k_j},x+x_j)-\frac{1}{\lambda_{k_j}^d}\, V\left(\frac{x}{\lambda_{k_j}} \right)\right|\dd x=0\;,
\end{equation*}
where $\lambda_k:=\|u(t_k)\|_m^{-m/(d-2)}$ and $V$ is the unique
radially symmetric minimiser of $\FF$ in $\YY_{M_c}$ such that
$\|V\|_m=1$. Assume further that
\begin{equation*}
\label{eq:bdm2} \mathcal{M}_2 := \sup_{t\in [0,T)} \int_{\RR^d}
|x|^2\, u(t,x)\, \dd x < \infty\;,
\end{equation*}
then
\begin{equation}
\label{eq:cvxkcom}
\lim_{j \to \infty} x_{j} = \bar{x} \;\;\mbox{ where }\;\; \bar{x} := \frac{1}{M_c}\, \int_{\RR^d} x\, u_0(x)\, \dd x\;.
\end{equation}
\end{proposition}
Since $\mu_{M_1+M_2}=\mu_{M_1}+\mu_{M_2}$ for $M_1\le M_c$ and  $M_2\le M_c$, the concentration compactness result as stated by P.-L. Lions~\cite{lions84} does not seem to apply directly. However, we follow the approach of M.~Weinstein~\cite{Weinstein86} to prove that the conclusion still holds true.
\medskip

\proof We set $v_k(x):=\lambda_k^d\, u(t_k,\lambda_k\,x)$ and aim at proving that $(v_k)_k$ converges strongly in $\LL^1(\RR^d)$. For this purpose we employ in Step~1 the concentration-compactness principle~\cite[Theorem~II.1]{lions84} to show that $(v_k)_k$ is tight up to translations. We argue in Step~2 as in~\cite[Theorem~1]{Weinstein86} to establish that $(v_k)_k$ has a limit in $\LL^1(\RR^d)$ and identify the limit. In the last step we use the additional bound on the second moment to show that the dynamics does not escape at infinity.

\noindent{\bf Step 1 - Tightness.} Obviously,
\begin{equation}
\label{ondine}
\|v_k\|_1 = M_c \;\;\mbox{ and }\;\; \|v_k\|_m = 1 \;\;\mbox{ for }\;\; k\ge 1\;.
\end{equation}
The concentration-compactness principle \cite{lions84} implies that there exists a sub-sequence (not relabelled) satisfying one of the three following properties:
\begin{description}
\item[(Compactness)] There exists a sequence $(a_k)_k$ in $\RR^d$ such that $(v_k(\cdot + a_k))_k \in \RR^d$ is tight, that is, for each $\eps>0$ there is $R_\eps>0$ such that
  \begin{equation}
    \label{toulouse}
    \int_{B(a_k,R_\eps)} v_k(x) \dd x \ge M_c - \eps\;.
  \end{equation}
\item[(Vanishing)] For all $R \ge 0$
  \begin{equation}
    \label{eq:5}
    \lim_{k \to \infty} \sup_{y \in \RR^d} \int_{B(y,R)}v_k(x)\dx =0\;.
  \end{equation}
\item[(Dichotomy)] There exists $\mu \in (0,M_c)$ such that for all
$\eps>0$, there exist $k_0\ge1$ and three sequences of non-negative,
integrable and compactly supported functions $(y_k^\eps)_k$,
$(z_k^\eps)_k$, and $(w_k^\eps)_k$ satisfying
$v_k=w_k^\eps+y_k^\eps+z_k^\eps$, 
  \begin{equation}\label{eq:6}
\left\{
    \begin{array}{l}
       \big|\|y_k^\eps\|_1-\mu\big|\le
\eps\;,\big|\|z_k^\eps\|_1-(M_c-\mu) \big| \le \eps\;,\|w_k^\eps\|_1
\le \eps\;,\vspace{.3cm}\\ 
\displaystyle \lim_{k \to \infty} \mbox{dist } (\mbox{supp }
y_k^\eps,\mbox{supp } z_k^\eps )=\infty\;, 
    \end{array}
\right.
  \end{equation}
for any $k \ge k_0$.
\end{description}

\medskip

As usual we shall rule out the possible occurrence of vanishing and
dichotomy. To this end we argue as
in~\cite[Theorem~II.1]{lions84}. Let us first notice that by the
scaling and non-negativity properties of the free energy, \eqref{eq:1}
and \eqref{eq:27}, $\FF[u(t_k)]\in \left[ 0,\FF[u_0] \right]$ and 
\begin{equation}\label{eq:ffvkto0}
 \lim_{k \to \infty} \FF[v_k]=\lim_{k \to
\infty}\|u(t_k)\|_m^{-m}\,\FF[u(t_k)] =0\;. 
\end{equation}
Consequently, since $\|v_k\|_m=1$ by the definition of $\lambda_k$, we have
\begin{equation}
  \label{eq:8}
  \lim_{k \to \infty} \WW(v_k) =  \lim_{k \to \infty}
\frac{2}{c_d}\left(\frac{1}{m-1} \|v_k\|_m^m -  \FF[v_k] \right) =
\frac{2}{c_d\,(m-1)}>0 \;. 
\end{equation}

\smallskip

\noindent$\bullet$ Let us first show that vanishing does not take
place and argue by contradiction. We split the non-local term
$\WW(v_k)$ in three parts. If $|x-y|$ is small, we control the
corresponding term by the bound in $\LL^1\cap\LL^m$ of $v_k$. If
$|x-y|$ is large the corresponding term is controlled by the
$\LL^1$-bound of $v_k$. And the remaining term converges to zero if we
assume that vanishing occurs which contradicts \eqref{eq:8}. Indeed,
if $q \in \left((d-1)/(d-2),d/(d-2) \right)$ and $R>0$, it follows
from the H\"older and Young inequalities that 
\begin{align*}
  \WW(v_k) &= \iint_{\RR^d\times\RR^d}
\frac{v_k(x)\,v_k(y)}{|x-y|^{d-2}} \un_{\left[
0,1/R\right]}(|x-y|)\dy\dx + \iint_{\RR^d\times\RR^d}
\frac{v_k(x)\,v_k(y)}{|x-y|^{d-2}} \un_{\left(1/R,R
\right)}(|x-y|)\dy\dx \\ 
&\hspace{1cm}+   \iint_{\RR^d\times\RR^d}
\frac{v_k(x)\,v_k(y)}{|x-y|^{d-2}}
\un_{\left[R,\infty\right)}(|x-y|)\dy\dx\\ 
&\le \|v_k\|^2_{2q/(2q-1)}\left( \int_{\RR^d} |x|^{-q(d-2)}
\un_{\left[ 0,{1}/{R}\right]}(|x|)\dx\right)^{1/q} +
R^{d-2}\int_{\RR^d} v_k(x)\int_{B(x,R)} v_k(y)\dy\dx \\ 
&\hspace{1cm} +   \frac{1}{R^{d-2}}\left( \int_{\RR^d}
v_k(x)\dx\right)^2\\ 
&\le
\|v_k\|_m^{m/[q(m-1)]}\,\|v_k\|_1^{[(2q-1)/q]-[d/(q(d-2))]}\,\left(
d\,\sigma_d\int_0^{1/R} r^{d-1-q(d-2)} \dr\right)^{1/q} \\ 
&\hspace{1cm} + R^{d-2}M_c \sup_{x \in \RR^d}\int_{B(x,R)}
v_k(y)\dy\dx + \frac{M_c^2}{R^{d-2}}\\ 
&\le C(q) \,\frac{1}{R^{[d-q(d-2)]/q}} + R^{d-2}M_c \sup_{x \in
\RR^d}\int_{B(x,R)} v_k(y)\dy\dx + \frac{M_c^2}{R^{d-2}} 
\end{align*}
We let $k\to\infty$ in the above inequality and use the vanishing
assumption~\eqref{eq:5} to obtain that 
\begin{equation*}
  \limsup_{k \to \infty} \WW(v_k) \le C(q) \left(
R^{2-d}+R^{-(d-q(d-2))/q}\right)\;. 
\end{equation*}
We next let $R$ to infinity to conclude that $\WW(v_k)$ converges to
zero as $k\to\infty$ which contradicts~\eqref{eq:8}. 
\smallskip

\noindent$\bullet$ Let us next assume for contradiction that dichotomy
takes place. We have 
\begin{equation*}
  \WW(v_k) - \WW(y_k^\eps) - \WW(z_k^\eps) =- \WW(w_k^\eps) + {\rm
I_1} + {\rm I_2}\;, 
\end{equation*}
where
\begin{equation*}
 {\rm I_1}:= 2\iint_{\RR^d\times\RR^d} y_k^\eps(x)\,z_k^\eps(y)\,|x-y|^{2-d}\dx\dy
\;\mbox{ and }\;
 {\rm I_2}:= 2\iint_{\RR^d\times\RR^d} v_k(x)\,w_k^\eps(y)\,|x-y|^{2-d}\dx\dy  \;.
\end{equation*}
On the one hand, setting $d_k^\eps:=\mbox{dist} (\mbox{supp } y_k^{\eps},\mbox{supp } z_k^{\eps})$, we have
\begin{align*}
 |{\rm I_1}| &\le \iint_{\RR^d \times \RR^d} y_k^\eps(x)\,z_k^\eps(y)\,\un_{(0,d_k^\eps)}(|x-y|)\,|x-y|^{2-d}\dx\dy\\
 &\qquad+ \iint_{\RR^d \times \RR^d} y_k^\eps(x)\,z_k^\eps(y)\,\un_{[d_k^\eps,\infty)}(|x-y|)\,|x-y|^{2-d}\dx\dy \;.
\end{align*}
Thanks to the definition of $d_k^\eps$ the first integral vanishes and we arrive at
\begin{equation*}
  |{\rm I_1}|\le M_c^2\,(d_k^\eps)^{2-d} \;.
\end{equation*}
On the other hand it follows from \eqref{ondine}, \eqref{eq:6}, the HLS inequality \eqref{eq:HLSineq} applied to $f=v_k$, $g=w_k^\eps$, $\lambda=d-2$ and $p=q=2d/(d+2)$, and the H\"older inequality with $1<2d/(d+2)<m$ that
\begin{align*}
  |{\rm I_2}| & \le C_{\rm HLS}\, \|v_k\|_{2d/(d+2)}\, \| w_k^\eps\|_{2d/(d+2)} \le C_{\rm HLS}\, \|v_k\|_m^{m/2}\, \|v_k\|_1^{1/d}\, \| w_k^\eps\|_m^{m/2} \, \| w_k^\eps\|_1^{1/d}  \\
  & \le C_{\rm HLS}\, M_c^{1/d}\,\|w_k^\eps\|_{m}^{m/2}\,\eps^{1/d}\;,
\end{align*}
and $0\le w_k^\eps \le v_k$ and \eqref{ondine} imply that $\|w_k^\eps\|_m\le 1$.
Similarly by the variant of the HLS inequality \eqref{eq:adaptedHLSineq}, we obtain
\begin{equation*}
  |\WW(w_k^\eps)|\le C_*\,\|w_k^\eps\|_m^{m}\,\eps^{2/d} \le C_*\, \eps^{2/d}\;.
\end{equation*}
Combining these estimates, we have thus shown that, given $\eps\in (0,1)$, there exists $k_\eps\ge \eps^{-1}$ such that
\begin{equation}
  \label{eq:10}
  |\WW(v_{k_\eps}) - \WW(y_{k_\eps}^\eps) -\WW(z_{k_\eps}^\eps)| \le \eps^{1/d}\;.
\end{equation}
Since $w_{k_\eps}^\eps$ is non-negative and the supports of $y_{k_\eps}^\eps$ and $z_{k_\eps}^\eps$ are disjoint we have
$$
\| y_{k_\eps}^\eps + z_{k_\eps}^\eps + w_{k_\eps}^\eps \|_m^m \ge \| y_{k_\eps}^\eps + z_{k_\eps}^\eps \|_m^m = \| y_{k_\eps}^\eps \|_m^m + \|  z_{k_\eps}^\eps \|_m^m\;,
$$
and we deduce from \eqref{eq:10} that
\begin{align*}
  \FF[v_{k_\eps}] &= \frac{1}{m-1} \|y_{k_\eps}^\eps + z_{k_\eps}^\eps + w_{k_\eps}^\eps\|_m^m - \frac{c_d}{2} \WW(v_{k_\eps})\\
& \ge \frac{1}{m-1} \left( \|y_{k_\eps}^\eps\|_m^m + \|z_{k_\eps}^\eps\|_m^m \right) - \frac{c_d}{2} \WW(y_{k_\eps}^\eps) - \frac{c_d}{2} \WW(z_{k_\eps}^\eps)- \frac{c_d}{2}\,\eps^{1/d}\\
& \ge \FF[y_{k_\eps}^\eps] + \FF[z_{k_\eps}^\eps] - \frac{c_d}{2}\, \eps^{1/d}\;.
\end{align*}
The above inequality, \eqref{eq:6}, \eqref{eq:ffvkto0}, and the non-negativity of $\FF$ for functions with $\LL^1$-norm lower or equal to $M_c$ then entail that
\begin{equation}
\label{romeo}
  \lim_{\eps \to 0} \FF[y_{k_\eps}^\eps]= \lim_{\eps \to 0} \FF[z_{k_\eps}^\eps] =0 \;.
\end{equation}
Now, \eqref{eq:lmbound} and \eqref{eq:6} imply
\begin{equation*}
  0 = \lim_{\eps \to 0}\FF[y_{k_\eps}^\eps] \ge \lim_{\eps \to 0} \frac{C_*\,c_d}{2} \left( M_c^{2/d} - \|y_{k_\eps}^\eps\|_1^{2/d} \right)\, \|y_{k_\eps}^\eps\|_m^m \ge \frac{C_*\,c_d}{2} \left( M_c^{2/d} - \mu^{2/d} \right)\, \lim_{\eps \to 0} \|y_{k_\eps}^\eps\|_m^m\;,
\end{equation*}
and a similar inequality for $z_{k_\eps}^\eps$ (with $M_c-\mu$ instead of $\mu$), hence
\begin{equation}
\label{juliette}
  \lim_{\eps \to 0} \|y_{k_\eps}^\eps\|_m^m  = \lim_{\eps \to 0} \|z_{k_\eps}^\eps\|_m^m =0 \;.
\end{equation}
Combining \eqref{romeo} and \eqref{juliette} gives
\begin{equation*}
  \lim_{\eps \to 0} \WW(y_{k_\eps}^\eps)=\lim_{\eps \to 0} \WW(z_{k_\eps}^\eps) =0\;,
\end{equation*}
which, together with~\eqref{eq:10}, implies that $\left( \WW(v_{k_\eps})\right)_\eps$ goes to 0 as $\eps$ goes to infinity and contradicts~\eqref{eq:8}.

\medskip

Having excluded the vanishing and dichotomy phenomena we thus conclude that there exists a sequence $(a_k)_k$ in $\RR^d$ such that $(v_{k}(\cdot + a_k))_k$ is tight, that is, satisfies \eqref{toulouse}.

\medskip

\noindent{\bf Step~2 - Compactness in $\LL^m$.} We now aim at showing that a sub-sequence of $(v_k(\cdot + a_k))_k$ converges in $\LL^1(\RR^d) \cap \LL^m(\RR^d)$ towards a minimiser of $\FF$ in $\YY_{M_c}$. We set $V_k(x):=v_k(x+a_k)$ for $x \in \RR^d$ and $k \ge 1$. By virtue of \eqref{ondine} we may assume (after possibly extracting a sub-sequence) that there is a non-negative $V_\infty \in \LL^m(\RR^d)$ such that
\begin{equation}
  \label{eq:12}
  V_k \rightharpoonup V_\infty \quad \mbox{weakly in $\LL^m(\RR^d)$} \;.
\end{equation}
By \eqref{ondine}, \eqref{toulouse}, and \eqref{eq:12} we have $V_\infty$ is non-negative with $\|V_\infty\|_1=M_c$ and $\|V_\infty\|_m \le 1$.

To prove the convergence of $\WW(V_k)$ to $\WW(V_\infty)$, we proceed as in Step~1 and split $\RR^d \times \RR^d$ in three parts. If $q \in \left((d-1)/(d-2),d/(d-2) \right)$ we have
\begin{align*}
  |\WW(V_k)-\WW(V_\infty)|
& \le \frac{2\,M_c^2}{R^{d-2}} + C(q)\, R^{[q(2-d)+d]/q} \\
&\hspace{1cm}+ \left| \iint_{\RR^d \times \RR^d} \left[ V_k(x)\,V_k(y)-V_\infty(x)\,V_\infty(y) \right]\, \frac{\un_{\left(1/R,R\right)}(|x-y|)}{|x-y|^{d-2}}\dy\dx \right|\;.
\end{align*}
Since $x \mapsto \un_{\left(1/R,R\right)}\,(|x|) |x|^{2-d} \in \LL^1(\RR^d) \cap \LL^\infty (\RR^d)$, the weak convergence \eqref{eq:12} ensures that $(x,y) \mapsto V_k(x)\,V_k(y)$ converges weakly toward $(x,y) \mapsto V_\infty(x)\,V_\infty(y)$ in $\LL^m(\RR^d \times \RR^d)$ so that the last term of the right-hand side converges to zero as $k\to \infty$. Therefore
\begin{equation*}
  \limsup_{k \to \infty} |\WW(V_k)-\WW(V_\infty)| \le C(q) \left( R^{2-d}+R^{-[d-q(d-2)]/q}\right)\;.
\end{equation*}
We then let $R\to\infty$ to obtain
\begin{equation*}
  \lim_{k \to \infty} \WW(V_k) =\WW(V_\infty)\;.
\end{equation*}

Owing to the lower semi-continuity of the $\LL^m$-norm and \eqref{eq:ffvkto0} we have
\begin{equation*}
  \FF[V_\infty]\le \frac{1}{m-1} \liminf_{k \to \infty} \|V_k\|_m^m - \frac{c_d}{2} \lim_{k \to \infty} \WW(V_k)\le \lim_{k \to \infty} \FF[V_k] = 0\;,
\end{equation*}
while Proposition~\ref{prop:inffreenrj} warrants that $\FF[V_\infty]\ge 0$ as $V_\infty\in\YY_{M_c}$. Consequently, $\FF[V_\infty]=0$ and the strong convergence of $(V_k)_k$ to $V_\infty$ in $\LL^m(\RR^d)$ readily follows: indeed,
\begin{align*}
  \frac{1}{m-1} \|V_\infty\|_m^m &= \FF[V_\infty] + \frac{c_d}{2}\WW(V_\infty) = \lim_{k \to \infty} \left(\FF[V_k] + \frac{c_d}{2}\WW(V_k) \right) = \frac{1}{m-1} \lim_{k \to \infty} \|V_k\|_m^m\;.
\end{align*}
We have thus shown that $V_\infty$ is a minimiser of $\FF$ in
$\YY_{M_c}$ with the additional property $\|V_\infty\|_m=1$.
Furthermore, according to the characterisation of
the minimisers given in
Proposition~\ref{prop:caracterisationminimiser}, there exists $y_0
\in \RR^d$ such that $V_\infty(\cdot + y_0)=:V$ is the unique
radially symmetric minimiser of $\FF$ in $\YY_{M_c}$ with
$\|V\|_m=1$. Coming back to the original variables we have proved
that $\left( x \mapsto \lambda_k^d\, u(t_k,\lambda_k\,
(x+a_k+y_0)) \right)_k$ converges to $V$ in $\LL^1(\RR^d)$ and
$\LL^m(\RR^d)$. Setting $x_k=\lambda_k\,(a_k+y_0)$ gives
\begin{equation}
\label{lille}
\lim_{k\to\infty} \int_{\RR^d} \left\vert u(t_k,x+x_k) - \frac{1}{\lambda_k^d}\, V\left( \frac{x}{\lambda_k}\right) \right\vert\, \dd x = 0\;,
\end{equation}
and thus the first assertion of Proposition~\ref{prop:how}.

\medskip
\noindent{\bf Step~3 - Convergence of $(x_k)_k$.} We first note that
$$
\int_{\RR^d} x\, u(t,x)\,\dd x = \int_{\RR^d} x\, u_0(x)\,\dd x
$$
for $t\in [0,T_\omega)$ so that we have also
$$
\bar{x} = \frac{1}{M_c}\, \int_{\RR^d} x\, u(t,x)\,\dd x
$$
for $t\in [0,T_\omega)$ by \eqref{evian}. Next, for $\eps\in (0,1)$, we have
\begin{align*}
\left\vert (\bar{x} - x_k)\, \int_{B(x_k,\eps)} u(t_k,x)\,\dd x \right\vert &\le \left\vert \int_{B(x_k,\eps)} (\bar{x} - x)\, u(t_k,x)\,\dd x \right\vert + \left\vert \int_{B(x_k,\eps)} (x - x_k)\, u(t_k,x)\,\dd x \right\vert \\
& \le \int_{\{|x-x_k|\ge\eps\}} \vert \bar{x} - x\vert\, u(t_k,x)\,\dd x + \eps\,M_c\\
& \le \int_{\{|x-x_k|\ge\eps\,,\,|x-\bar{x}|\le 1/\eps\}} \vert \bar{x} - x\vert\, u(t_k,x)\,\dd x \\
& + \int_{\{|x-x_k|\ge\eps\,,\,|x-\bar{x}|> 1/\eps\}} \vert \bar{x} - x\vert\, u(t_k,x)\,\dd x + \eps\,M_c\\
& \le \frac{1}{\eps}\, \int_{\{|x-x_k|\ge\eps\}} u(t_k,x)\,\dd x + \eps\, (\mathcal{M}_2 + M_c)\\
& \le \frac{1}{\eps}\, \int_{\RR^d} \left\vert u(t_k,y+x_k) - \frac{1}{\lambda_k^d}\, V\left( \frac{y}{\lambda_k} \right) \right\vert\, \dd y \\
& + \frac{1}{\eps}\, \int_{\{|x|\ge \eps/\lambda_k\}} V(y)\dd y + \eps\, (\mathcal{M}_2 + M_c)\;.
\end{align*}
Since $\lambda_k\to 0$ as $k\to\infty$ we infer from \eqref{lille} and the integrability of $V$ that
$$
\limsup_{k\to\infty} \left\vert (\bar{x} - x_k)\, \int_{B(x_k,\eps)} u(t_k,x)\,\dd x \right\vert \le \eps\, (\mathcal{M}_2 + M_c)\;.
$$
Using once more \eqref{lille} we readily deduce that
$$
\lim_{k\to\infty} \int_{B(x_k,\eps)} u(t_k,x)\,\dd x = \int_{\RR^d} V(x) \,\dd x = M_c\;.
$$
Combining the previous two limits gives
$$
M_c\, \limsup_{k\to\infty} \left\vert \bar{x} - x_k \right\vert \le \eps\, (\mathcal{M}_2 + M_c)\;,
$$
whence the last assertion of Proposition~\ref{prop:how} by letting $\eps\to 0$. \finproof
For radially symmetric solutions we can remove the additional assumption on the second moment.
\begin{corollary}[Radially symmetric blow-up]
Let $u_0$ be a radially symmetric initial condition satisfying \eqref{Eqn:Assumptions} with $\|u_0\|_1=M_c$ and consider a radially symmetric free energy solution $u$ to~\eqref{eq:spf} on $[0,T_\omega)$ with initial condition $u_0$ and $T_\omega \in (0,\infty]$ and such that $\|u(t)\|_m\to\infty$ as $t\nearrow T_\omega$. If $(t_k)_k$ is a sequence of positive real numbers such that $t_k\to T_\omega$ as $k\to\infty$, there is a sub-sequence $(t_{k_j})_j$ of $(t_k)_k$ such that
\begin{equation*}
  \lim_{j \to \infty} \int_{\RR^d}\left| u(t_{k_j},x)-\frac{1}{\lambda_{k_j}^d}\, V\left(\frac{x}{\lambda_{k_j}} \right)\right|\dd x=0\;,
\end{equation*}
where $\lambda_k:=\|u(t_k)\|_m^{-m/(d-2)}$ and $V$ is the unique radially symmetric minimiser of $\FF$ in $\YY_{M_c}$ such that  $\|V\|_m=1$.
\end{corollary}
\proof
The only modification of the proof of Proposition~\ref{prop:how} is to show that we can choose $a_k=0$ for all $k$ at the end of Step~1. Indeed, we claim that if $\eps \in (0,M_c/4)$ we have $|a_k|\le R_\eps$, where $a_k$ and $R_\eps$ are defined in~\eqref{toulouse}. Otherwise $B(a_k,R_\eps)$ and $B(-a_k,R_\eps)$ are disjoint and the radial symmetry of $v_k$ and~\eqref{toulouse} imply that
\begin{equation*}
  M_c \ge \int_{B(a_k,R_\eps) \cup B(-a_k,R_\eps)}v_k(x)\dx = 2 \int_{B(a_k,R_\eps)}v_k(x)\dx \ge 2\,(M_c - \eps) \ge \frac{3\,M_c}{2}
\end{equation*}
and a contradiction. Therefore $B(a_k,R_\eps) \subset B(0,2\,R_\eps)$ and thus
\begin{equation*}
   \int_{ B(0,2\,R_\eps)}v_k(x)\dx \ge M_c - \eps
\end{equation*}
by~\eqref{toulouse}.
\finproof
\subsubsection{When would it blow-up?}
\begin{proposition}[Global existence in the critical case]\label{prop:gloexistcritic}
Let $u_0$ be an initial condition satisfying~\eqref{Eqn:Assumptions} with $\|u_0\|_1=M_c$ and consider a free energy solution $u$ to \eqref{eq:spf} on $[0,T_\omega)$ with initial condition $u_0$ and $T_\omega \in (0,\infty]$. Then $T_\omega=\infty$.
\end{proposition}
The proof of this proposition relies on Proposition~\ref{prop:how} and the following control of the behaviour of free energy solutions for large $x$:
\begin{lemma}[Control of the tail]\label{lem:controltail}
Consider a free energy solution $u$ to \eqref{eq:spf} on $[0,T_\omega)$ with initial condition $u_0$ satisfying \eqref{Eqn:Assumptions} and $T_\omega\in (0,\infty]$. If $t\longmapsto \FF[u(t)]$ is bounded from below in $[0,T)$ for some $T \le T_\omega$ then
\begin{equation*}
  \lim_{R \to \infty} \; \sup_{t \in [0,T)} \int_{\{|x|>R\}} |x|^2\, u(t,x) \dd x=0\;.
\end{equation*}
\end{lemma}
\proof
Consider a non-decreasing function $\xi\in\mathcal{C}^\infty(\RR)$ such that $\xi(r)=0$ for $|r|\le 1$ and  $\xi(r)=1$ for $|r|\ge 2$ and define
\begin{equation*}
  \Phi_R(r)=r\,\xi^4\left(\frac{r}{R}\right) \;\;\mbox{ for }\;\; r\in\RR \;\;\mbox{ and }\;\; R>0\;.
\end{equation*}
The support of $\Phi_R$ is included in $\RR^d \setminus B(0,R)$ and, introducing
\begin{equation*}
  \II_R(t):= \int_{\{|x|>R\}} \left| \left(\frac{2\,m}{2\,m-1}\, \nabla u^{(2m-1)/2} - u^{1/2}\,\nabla\phi\right)(t,x) \right|^2\dx\;,
\end{equation*}
we have
\begin{align*}
  \frac{\dd}{\dt}\int_{\RR^d}  &\Phi_R\left(|x|^2\right)\, u(t,x)\dd x
  \\&=- \int_{\RR^d}2\,x \,\Phi_R'\left(|x|^2\right)\, \left(\frac{2\,m}{2\,m-1}\, \nabla u^{(2m-1)/2}-u^{1/2}\,\nabla \phi \right)(t,x)\dd x\\
&\le 2\left(\int_{\RR^d}|x|^2 \left|\Phi_R'\left(|x|^2\right)\right|^2\, u(t,x) \dd x\right)^{1/2}\, \II_R(t)^{1/2}\;.
\end{align*}
By the definition of $\Phi_R$, we have
\begin{align*}
   |\Phi_R'(r)|^2 &\le \left| \xi^4\left(\frac{r}{R}\right)+4\, \frac{r}{R}\, \xi^3\left(\frac{r}{R}\right)\, \xi'\left(\frac{r}{R}\right) \right|^2 \le  2\, \xi^8\left(\frac{r}{R}\right) + 32\, \left(\frac{r}{R}\right)^2\, \xi^6\left(\frac{r}{R}\right)\, \left|\xi'\left(\frac{r}{R}\right)\right|^2\\
&\le  \xi^4\left(\frac{r}{R}\right)\, \left( 2+32\,\sup_{z\in\RR} \left|z\,\xi'\left(z\right)\right|^2\right)\;,
\end{align*}
so that $r\,|\Phi_R'(r)|^2 \le C\,\Phi_R(r)$ for $r\in\RR$. Therefore, for any $t \in [0,T)$,
$$
\frac{\dd}{\dt} \int_{\RR^d} \Phi_R\left(|x|^2\right)\, u(t,x)\dd x \le C\, \left( \int_{\RR^d} \Phi_R\left(|x|^2\right)\, u(t,x)\dd x \right)^{1/2}\, \II_R(t)^{1/2}\;,
$$
hence
\begin{equation*}
 \left(\int_{\RR^d} \Phi_R\left(|x|^2\right)\, u(t,x) \dx\right)^{1/2} \le \left(\int_{\RR^d} \Phi_R\left(|x|^2\right)\, u_0(x) \dx \right)^{1/2} + \frac{C}{2}\, \int_0^t \II_R^{1/2}(s)\ds\;.
\end{equation*}
Now, since $\FF[u(t)]$ is bounded from below in $[0,T)$, we have
\begin{equation*}
\int_0^T \int_{\RR^d}  \left|\left(\frac{2\,m}{2\,m-1}\, \nabla u^{(2m-1)/2} - u^{1/2}\,\nabla \phi\right)(s,x) \right|^2 \dx \ds \le \sup_{t\in [0,T)}{\{ \FF[u_0]-\FF[u(t)] \}} < \infty\;,
\end{equation*}
so that
\begin{equation*}
  \lim_{R \to \infty} \int_0^T \II_R^{1/2}(s)\ds =0
\end{equation*}
by the Lebesgue dominated convergence theorem. Therefore,
\begin{equation*}
  \limsup_{R \to \infty} \sup_{t \in [0,T)} \int_{\RR^d} \Phi_R\left(|x|^2\right)\, u(t,x) \dx = 0\,,
\end{equation*}
from which the lemma follows.
\finproof

\noindent {\sl Proof of Proposition~{\rm \ref{prop:gloexistcritic}}.\/} Assume for contradiction that $T_\omega$ is finite and let $(t_k)_k$ be a sequence of positive real numbers such that $t_k\to T_\omega$ as $k\to\infty$. Observe that Theorem~\ref{charac} entails that $\|u(t)\|_m\to\infty$ as $t\to T_\omega$. On the one hand we infer from the nature of the blow-up given in Proposition~\ref{prop:how} that there are a sub-sequence of $(t_k)_k$ (not relabelled) and a sequence $(x_k)_k$ in $\RR^d$ such that
\begin{eqnarray}
\label{eq:acdc1}
& & \lim_{k\to\infty} x_k = \bar{x}:= \frac{1}{M_c}\, \int_{\RR^d} x\, u_0(x)\, \dd x\;,\\
\label{eq:acdc2}
& &   \lim_{k \to \infty} \int_{\RR^d}\left| u(t_k,x+x_k)-\frac{1}{\lambda_k^d}\, V\left(\frac{x}{\lambda_k} \right)\right|\dd x=0
\end{eqnarray}
with $\lambda_k:=\|u(t_k)\|_m^{-m/(d-2)}$. On the other hand it follows from Proposition~\ref{prop:inffreenrj} and Lemma~\ref{lem:viriel} that $\FF[u(t)]\ge 0$ for $t\in [0,T_\omega)$ so that
\begin{eqnarray}
\nonumber
\int_{\RR^d} |x-\bar{x}|^2\, u(t,x)\,\dd x & = & \int_{\RR^d} |x-\bar{x}|^2\, u_0(x)\,\dd x + 2(d-2)\, \int_0^t \FF[u(s)]\,\dd s \\
\label{eq:acdc3}
& \ge & \int_{\RR^d} |x-\bar{x}|^2\, u_0(x)\,\dd x > 0 \;,
\end{eqnarray}
and Lemma~\ref{lem:controltail} may be applied to obtain
\begin{equation}
\label{eq:acdc4}
\lim_{R \to \infty} \; \sup_{t \in [0,T_\omega)} \int_{\{|x|>R\}} |x|^2\,u(t,x)\dd x=0\;.
\end{equation}

Now, for $k\ge 1$ and $R\ge |\bar{x}|$ we have
\begin{align*}
\int_{\RR^d} |x-\bar{x}|^2\, u(t_k,x)\,\dd x & \le 2\, \int_{\{|x-\bar{x}|\ge 2R\}} \left( |x|^2 + |\bar{x}|^2 \right)\, u(t_k,x)\,\dd x \\
& + \int_{\{|x-\bar{x}|<2R\}} |x-\bar{x}|^2\, \left[ u(t,x) - \frac{1}{\lambda_k^d}\, V\left( \frac{x-x_k}{\lambda_k} \right) \right]\, \dd x \\
& + \int_{\{|x-\bar{x}|<2R\}} \frac{|x-\bar{x}|^2}{\lambda_k^d}\, V\left( \frac{x-x_k}{\lambda_k} \right)\, \dd x \\
& \le 4\, \int_{\{|x|\ge R\}} |x|^2\, u(t_k,x)\, \dd x + R^2\, \int_{\RR^d} \left\vert u(t,x) - \frac{1}{\lambda_k^d}\, V\left( \frac{x-x_k}{\lambda_k} \right) \right\vert\, \dd x \\
& + \int_{\RR^d} \left\vert \lambda_k\, x + x_k - \bar{x} \right\vert^2\, V(x)\, \dd x \\
& \le 4\, \sup_{t\in [0,T_\omega)} \int_{\{|x|\ge R\}} |x|^2\, u(t,x)\, \dd x + R^2\, \int_{\RR^d} \left\vert u(t,x) - \frac{1}{\lambda_k^d}\, V\left( \frac{x-x_k}{\lambda_k} \right) \right\vert\, \dd x \\
& + 2\, \left\vert x_k - \bar{x} \right\vert^2\, M_c + 2\,\lambda_k^2\, \int_{\RR^d} |x|^2\, V(x)\, \dd x\;.
\end{align*}
Owing to \eqref{eq:acdc1}, \eqref{eq:acdc2}, and the convergence of $(\lambda_k)_k$ to zero we may let $k\to\infty$ in the previous inequality to obtain
$$
\limsup_{k\to\infty} \int_{\RR^d} |x-\bar{x}|^2\, u(t_k,x)\,\dd x \le 4\, \sup_{t\in [0,T_\omega)} \int_{\{|x|\ge R\}} |x|^2\, u(t,x)\, \dd x\;.
$$
We next pass to the limit as $R\to\infty$ with the help of \eqref{eq:acdc4} to conclude that
$$
\lim_{k\to\infty} \int_{\RR^d} |x-\bar{x}|^2\, u(t_k,x)\,\dd x = 0\;,
$$
which contradicts \eqref{eq:acdc3}. \finproof
\subsubsection{Does it blow-up?}\label{sec:does}
Let us first note that Proposition~\ref{prop:how} allows us to
describe the nature of the blow-up when it occurs. We define the
two following statements: 
\begin{equation}
  \label{p1}
 \mbox{There exists $(t_k)_k\nearrow\infty$ such that $(\|u(t_k)\|_m)_k$ is bounded} \tag{S1}
\end{equation}
\begin{equation}
  \label{p2}
 {\cal M}_2^\infty:=\lim_{t \to \infty} \int_{\RR^d}|x|^2 u(t,x)\dx < \infty \tag{S2}
\end{equation}

\begin{itemize}
\item If [not (S1)] and (S2): By Proposition~\ref{prop:how}, the
solution blows up as a Dirac mass at the centre of mass as $t$
goes to infinity. Moreover, the blow-up profile is described by
the minimisers of $\FF$ for the critical mass.

\item If (S1) and (S2): By the virial identity Lemma
\ref{lem:viriel}, $\FF[u(t_k)] \to 0$ so that $(u(t_k))_k$ is a
minimising sequence for $\FF$ in $\YY_{M_c}$. We expect that it
converges to the minimiser of $\FF$ in $\YY_{M_c}$ with centre of
mass $\bar x$ defined in~\eqref{eq:cvxkcom} and second moment
${\cal M}_2^\infty$.

\item If [not (S1)] and [not (S2)]: By Proposition~\ref{prop:how},
the solution blows up as a Dirac mass. However, we cannot prevent
the escape at infinity of the Dirac mass.

\item If (S1) and [not (S2)]: No precise information can be
deduced in this case. We cannot even rule out the possibility of
the existence of another sequence of times for which the
$\LL^m$-norm diverges.
\end{itemize}

In the radially symmetric case, if the initial condition is less concentrated than one of the stationary solutions, then we strongly believe that such a property remains true for all times, thus excluding the formation of a Dirac mass. According to the above discussion this prevents the blow-up of the $\LL^m$-norm in infinite time and give an example where (S1) and (S2) hold true. This is in sharp contrast with the two-dimensional PKS case where infinite time blow-up always occurs, see~\cite{bkln,BCM}.
\section{Sub-critical self-similar solutions}\setcounter{equation}{0}
The aim of this section is to prove the existence of self-similar solution by variational techniques. Actually, it is equivalent to show the existence of minimisers for the free energy $\GG$ associated to the rescaled problem~\eqref{eq:spscaled} given by
\begin{equation*}
  \GG[h] := \FF[h]+\frac{1}{2}M_2[h] \;\;\mbox{ with }\;\; M_2[h]:=\int_{\RR^d} |x|^2\, |h(x)| \dd x
\end{equation*}
for $h \in \LL^1(\RR^d; (1+x^2)\dd x) \cap \LL^m(\RR^d)$.
For $M>0$, we define
\begin{equation*}
  \nu_M:=\inf\{\GG[h]\,:\,h \in \ZZ_M\}\quad \mbox{with}\quad \ZZ_M:=\{h \in \YY_M\,:\, M_2[h]<\infty\}\;.
\end{equation*}

We first establish the following analogue of~Proposition~\ref{prop:inffreenrj}.
\begin{proposition}[Infimum of the rescaled free energy]
For $M>0$ and $h\in\ZZ_M$ we have
\begin{equation}
\label{eq:9}
\GG[h] \ge \frac{C_*\,c_d}{2}\, \left( M_c^{2/d} -M^{2/d} \right )\, \|h\|_m^m + \frac{1}{2} M_2[h]\;.
\end{equation}
In addition,
\begin{equation*}
  \left\{
    \begin{array}{ll}
     \nu_{M} >0 &\quad \mbox{if $M<M_c$}\;,\vspace{.2cm}\\
      \nu_{M_c} =0\;,&\vspace{.2cm}\\
\nu_{M} = -\infty &\quad \mbox{if $M>M_c$}\;.
    \end{array}
\right.
\end{equation*}
\end{proposition}
\proof The inequality \eqref{eq:9} readily follows from \eqref{eq:lmbound} and the definition of $\GG$. Consider next $M\ge M_c$ and put
\begin{equation*}
  h_R(x):=
\left\{
  \begin{array}{ll}
\displaystyle    \frac{M}{M_c} \frac1{R^d} \zeta^{d/(d-2)}\left(\frac{x}{R^d} \right)\quad &\mbox{if $x \in B(0,R)$}\;, \vspace{.3cm}\\
0\quad &\mbox{if $x \in \RR^d \setminus B(0,R)$}\;,
  \end{array}
\right.
\end{equation*}
where the function $\zeta$ is defined in Proposition~\ref{prop:caracterisationminimiser} and $R>0$. We compute $\GG[h_R]$ and use the property $\FF[\zeta^{d/(d-2)}]=0$ to obtain
\begin{equation*}
  \nu_M \le \GG[h_R]= \left(\frac{M}{M_c} \right)^2 R^{2-d} \left[ \frac{M_c}{2M} R^d M_2\left[ \zeta^{d/(d-2)} \right] - \left(1-\left(\frac{M_c}{M} \right)^{2-m} \right) \frac{\|\zeta^{d/(d-2)}\|_m^m}{m-1}\right]\;.
\end{equation*}
Now, either $M>M_c$ and the right-hand side of the above inequality diverges to $-\infty$ as $R\to 0$ since $d>2$ and $m<2$. Consequently $\nu_m=-\infty$ in that case. Or $M=M_c$ and we may let
$R\to 0$ in the above inequality to obtain that $\nu_{M_c}
\le 0$. Since $\GG$ is non-negative by Proposition~\ref{prop:inffreenrj}, we conclude that $\nu_{M_c}=0$.

Finally, assume for contradiction that $\nu_{M}=0$ for some $M <
M_c$ and let $(h_k)_k$ be a minimising sequence for $\GG$ in
$\ZZ_M$. Since $\GG[h_k] \ge \GG[|h_k|]$, $(|h_k|)_k$ is also a
minimising sequence for $\GG$ in $\ZZ_M$ and we infer
from \eqref{eq:9} that
\begin{equation*}
  \lim_{k \to \infty}\left(\|h_k\|_m + M_2[h_k] \right)=0\;.
\end{equation*}
By Vitali's theorem $(|h_k|)_k$ converges towards zero in
$\LL^1(\RR^d)$ which contradicts the fact that
$\|h_k\|_{\LL^1}=M$ for all $k \ge 1$. Therefore $\nu_{M} \neq 0$
and the non-negativity of $\GG$ in $\ZZ_M$ entails that
$\nu_M>0$. \finproof

\medskip

We next identify the minimisers of $\GG$ in $\ZZ_M$ for $M\in (0,M_c)$.
\begin{theorem}[Identification of minimisers]\label{thm:rescaledminimiser}
If $M\in (0,M_c)$ there is a unique minimiser $W_M$ of $\GG$ in
$\ZZ_M$. In addition, $W_M$ is non-negative radially symmetric and
non-increasing and there is a unique $\varrho_M>0$ such that
$W_M(x)=0$ for $|x|\ge\varrho_M$ and $\xi_M:=W_M^{m-1}$ solves 
$$
\Delta \xi_M + \frac{m-1}{m}\, \left( \xi_M^{1/(m-1)} + d \right) = 0
\;\;\mbox{ in }\;\; B(0,\varrho_M) \;\;\mbox{ with }\;\; \xi_M=0
\;\;\mbox{ on }\;\; \partial B(0,\varrho_M)\;. 
$$
\end{theorem}
Several steps are required to perform the proof of
Theorem~\ref{thm:rescaledminimiser} which borrows several arguments
from \cite{lieb83,liebyau}. We first establish the existence of
minimisers of $\GG$ in $\ZZ_M$ for $M\in (0,M_c)$. 
\begin{lemma}[Existence of minimisers]\label{lem:fssv32}
Consider $M\in (0,M_c)$. The functional $\GG$ has at least a minimiser
in $\ZZ_M$. In addition, every minimiser of $\GG$ in $\ZZ_M$ is
non-negative radially symmetric and non-increasing. 
\end{lemma}
\proof We first recall that, if $h\in \LL^1(\RR^d;(1+|x|^2)\dd x)$ and
$h^*$ denotes its symmetric decreasing rearrangement, then
$M_2[h^*]\le M_2[h]$. Thanks to this property, we may next argue as in
the proof of Lemma~\ref{lem:extremalfct} to conclude that there is at
least a minimiser of $\GG$ in $\ZZ_M$. 

Next, let $W$ be a minimiser of $\GG$ in $\ZZ_M$ and denote by $W^*$
its symmetric decreasing rearrangement. As 
$$
\|W^*\|_1=\| W\|_1\,, \;\; \|W^*\|_m=\| W\|_m \,, \;\;\mbox{ and }\;\;
M_2[W^*]\le M_2[W]\;, 
$$
$W^*$ belongs to $\ZZ_M$. In addition, by Riesz's rearrangement
inequality~\cite[Lemma~2.1]{lieb83},
$\WW(W)\le\WW(W^*)$. Consequently, $\nu_M=\GG[W]\ge \GG[W^*]$ and
$W^*$ is also a minimiser of $\GG$ in $\ZZ_M$. This last property
entails that 
$$
M_2[W^*]=M_2[W] \;\;\mbox{ and }\;\; \WW(W^*)=\WW(W)\;.
$$
Using once more \cite[Lemma~2.1~(ii)]{lieb83} we deduce from
$\WW(W^*)=\WW(W)$ that there is $y\in\RR^d$ such that $W(x)=W^*(x+y)$
for $x\in\RR^d$. Then $M_2[W^*]=M_2[W]$ implies that $y=0$, which
completes the proof. 
\finproof

\medskip

We are thus left with the uniqueness issue to complete the proof of
Theorem~\ref{thm:rescaledminimiser}. To this end we adapt the proof in
\cite[Section~IV.B]{liebyau} and first proceed as in the proof of
Proposition~\ref{prop:caracterisationminimiser} to identify the
Euler-Lagrange equation satisfied by the minimisers of $\GG$ in
$\ZZ_M$. 
\begin{lemma}\label{lem:fssv34}
Consider $M\in (0,M_c)$ and let $W$ be a minimiser of $\GG$ in
$\ZZ_M$. Then there is $\varrho>0$ such that $W(x)=0$ if
$|x|\ge\varrho$ and $\xi:=W^{m-1}$ is a non-negative radially
symmetric and non-increasing classical solution to 
$$
\Delta \xi + \frac{m-1}{m}\, \left( \xi^{1/(m-1)} + d \right) = 0
\;\;\mbox{ in }\;\; B(0,\varrho) \;\;\mbox{ with }\;\; \xi=0
\;\;\mbox{ on }\;\; \partial B(0,\varrho)\;. 
$$
In addition,
\begin{equation}
\label{eq:fssv310}
\frac{m}{m-1}\, W^{m-1} = \left( \EE \ast W - \frac{|x|^2}{2} +
\frac{1}{2} + \frac{m}{m-1}\, M^{m-1} - \frac{c_d}{M}\,\WW(W)
\right)_+\qquad\mbox{a.e. in $\RR^d$}\;. 
\end{equation}
\end{lemma}
Additional properties of minimisers of $\GG$ in $\ZZ_M$ can be deduced
from Lemma~\ref{lem:fssv34}. 
\begin{lemma}\label{lem:fssv34b}
Consider $M\in (0,M_c)$ and let $W$ be a minimiser of $\GG$ in $\ZZ_M$. Then
\begin{eqnarray}
\label{eq:fssvA2}
M_2[W] = (d-2)\,\FF[W] & = & 2(m-1)\,\nu_M\;, \\
\label{eq:fssvA3}
\frac{2m}{m-1}\, \|W\|_m^m + M_2[W] & = & \frac{2m}{m-1}\, M^m + M\;.
\end{eqnarray}
\end{lemma}
\proof We proceed as in \cite[Lemma~6]{liebyau}. By
Lemma~\ref{lem:fssv34} we have 
$$
-\frac{\dd}{\dd r} \left( r^{d-1}\, \frac{\dd\xi}{\dd r}(r) \right) =
\frac{m-1}{m}\, \left( r^{d-1}\, W(r) + d\, r^{d-1} \right) \;\;\mbox{
for }\;\; r\in (0,\varrho)\;, 
$$
where $\varrho$ denotes the radius of the support of $W$ and
$\xi:=W^{m-1}$. Introducing 
$$
Q(r) := \int_{B(0,r)} W(x)\,\dd x = \sigma_d\, \int_0^r W(z)\,
z^{d-1}\,\dd z \;\;\mbox{ for }\;\; r\in (0,\varrho)\;, 
$$
we integrate the previous differential equation to obtain
$$
-m\,r^{d-1}\,W(r)^{m-2}\,\frac{\dd W}{\dd r}(r) =
\frac{Q(r)}{\sigma_d} + r^d \;\;\mbox{ for }\;\; r\in (0,\varrho)\;. 
$$
Multiplying the above identity by $\sigma_d\, r\, W(r)$ and integrating
over $(0,\infty)$ then lead us to the formula 
$$
d\,\|W\|_m^m = \int_0^\infty r\,Q(r)\,W(r)\, \dd r + M_2[W]\,.
$$
As
$$
2\sigma_d\,\int_0^\infty r\,Q(r)\,W(r)\, \dd r = \WW(W)
$$
by Newton's theorem \cite[Theorem~9.7]{liebloss}, we end up with the
identity $(d-2)\,\FF[W]=M_2[W]$ and~\eqref{eq:fssvA2} follows by the
definition of $\nu_M$ and $\GG$. We next multiply \eqref{eq:fssv310}
by $2\,W$ and integrate over $\RR^d$ to obtain
\eqref{eq:fssvA3}. \finproof 

\medskip

We next prove the following comparison result.

\begin{lemma}\label{lem:fssv35}
Consider $M_1\in (0,M_c)$ and $M_2\in (0,M_c)$. For $i=1,2$ let $W_i$ be a minimiser of $\GG$ in $\ZZ_{M_i}$ and denote by $\varrho_i$ the radius of its support (which is finite according to Lemma~{\rm \ref{lem:fssv34}}). If $W_1(0) > W_2(0)$ then $Q_1(r) > Q_2(r)$ for $r\in (0,\max{\{\varrho_1,\varrho_2\}})$ where
$$
Q_i(r) := \int_{B(0,r)} W_i(x)\,\dd x \;\;\mbox{ for }\;\; r\in (0,\max{\{\varrho_1,\varrho_2\}}) \;\;\mbox{ and }\;\; i=1,2\;.
$$
\end{lemma}

Owing to Lemma~\ref{lem:fssv32} and Lemma~\ref{lem:fssv34}, the proof of Lemma~\ref{lem:fssv35} is similar to that of \cite[Lemma~10]{liebyau} to which we refer.

\medskip

\noindent {\sl Proof of Theorem~{\rm \ref{thm:rescaledminimiser}}.\/} Consider $M\in (0,M_c)$ and assume for contradiction that $\GG$ has two minimisers $W_1$ and $W_2$ in $\ZZ_M$ with $W_1(0)>W_2(0)$. Denoting by $\varrho_i$ the radius of the support of $W_i$ and introducing
$$
Q_i(r) := \int_{B(0,r)} W_i(x)\,\dd x
$$
for $r\in [0,\max{\{\varrho_1,\varrho_2\}}]$ and $i=1,2$, we infer from Lemma~\ref{lem:fssv35} that $Q_1(r)>Q_2(r)$ for all $r\in (0,\max{\{\varrho_1,\varrho_2\}})$. Then $\varrho_1\le\varrho_2$ and \eqref{eq:fssvA2} warrants that
\begin{align*}
2\,(m-1)\,\nu_M &= \sigma_d\, \int_0^\infty r^2\, \frac{\dd}{\dd r}\left( Q_i - M \right)(r)\,\dd r = \int_0^\infty 2\,r\, \left( M - Q_i(r) \right)\,\dd r
\end{align*}
for $i=1,2$. Consequently,
$$
\int_0^{\varrho_2} 2\,r\, \left( Q_1 - Q_2 \right)(r)\,\dd r = 0\;,
$$
which implies that $\varrho_1=\varrho_2$ and $Q_1=Q_2$, hence a contradiction. \finproof
\begin{corollary}\label{cor:last}
 If $M\in (0,M_c)$ there exists a self-similar solution $U_{M}$ to~\eqref{eq:spf} given by
 \begin{equation*}
   U_{M}(t,x) = \frac1{1+dt} W_M\left(\frac{x}{(1+dt)^{1/d}}\right)\;,
 \end{equation*}
 where $W_M$ is the unique minimiser of $\GG$ in $\ZZ_M$ given in Theorem~{\rm \ref{thm:rescaledminimiser}}.
\end{corollary}
\begin{remark}
Given $M \in (0,M_c)$, we expect that this self-similar solution
attracts the dynamics of~\eqref{eq:spf} for large times. Although
we can prove that the $\omega$-limit set of the rescaled
equation~\eqref{eq:spscaled} consists of stationary solutions, we
are yet lacking a uniqueness result to identify them as $W_{M}$.
\end{remark}
\noindent {\bf Acknowledgements.-} The authors are grateful to
Mohammed Lemou for pointing out~\cite{liebyau} and to Pierre Rapha\"el
for stimulating discussions. AB acknowledges
the support of bourse Lavoisier. JAC acknowledges the support from
DGI-MEC (Spain) project MTM2005-08024. AB and JAC acknowledge
partial support of the Acc. Integ./Picasso program HF2006-0198.
We thank the Centre de Recerca Matem\`atica
(Barcelona) for partial funding and for providing an excellent
atmosphere for research.


\bigskip\noindent{\small This paper is under the Creative Commons
licence Attribution-NonCommercial-ShareAlike 2.5.}
\end{document}